\documentclass[a4paper,11pt]{article}

\usepackage{amsmath,amsthm,amssymb,amsfonts}
\usepackage{mathtools} 
\usepackage{bm}
\usepackage{mathrsfs} 
\usepackage{changes} 
\usepackage{url} 
\usepackage[toc,page]{appendix} 
\usepackage{enumerate}
\usepackage{enumitem}

\usepackage{graphicx} 
\usepackage{longtable} 
\usepackage[all,knot,arc,import,poly]{xy}
\usepackage{newfloat} 
\DeclareFloatingEnvironment[fileext=frm,placement={ht},name=Frame]{algfloat}
\usepackage{caption}
\usepackage{placeins}
\usepackage{tikz}
\usepackage{subfig}
\usepackage{pgfplots}
\usepackage{color}
\usepackage{xcolor} 
\usepackage{array}
\usepackage{booktabs}
\usepackage{grffile}
\usetikzlibrary{patterns}
\usetikzlibrary{calc,trees,positioning,arrows,chains,shapes.geometric,decorations.pathreplacing,decorations.pathmorphing,shapes,matrix,shapes.symbols}
\usetikzlibrary{shapes,backgrounds}
\pgfplotsset{width=13cm,compat=1.14}

\usepackage{scalefnt}

\newcommand{\R}{\mathbb{R}}

\newcommand{\inner}[2]{\langle #1\,|\,#2\rangle} 
\newcommand{\Inner}[2]{\left\langle #1\,|\,#2\right\rangle}
\newcommand{\norm}[1]{\|{#1}\|}
\newcommand{\normq}[1]{\|{#1}\|^2}

\newcommand{\HH}{{\mathcal{H}}} 

\newcommand{\comenta}[1]{} 

\newcommand{\set}[1]{\{#1\}}

\newcommand{\lab}[1]{\label{#1}}

\newcommand{\iteref}[1]{\emph{\ref{#1}}}

\newcommand{\bprop}{\begin{proposition}}
\newcommand{\eprop}{\end{proposition}}
\newcommand{\blemm}{\begin{lemma}}
\newcommand{\elemm}{\end{lemma}}
\newcommand{\bdefi}{\begin{definition}}
\newcommand{\edefi}{\end{definition}}
\newcommand{\btheo}{\begin{theorem}}
\newcommand{\etheo}{\end{theorem}}
\newcommand{\bproo}{\begin{proof}}
\newcommand{\eproo}{\end{proof}}
\newcommand{\brema}{\begin{remark}}
\newcommand{\erema}{\end{remark}}
\newcommand{\bitem}{\begin{itemize}}
\newcommand{\eitem}{\end{itemize}}
\newcommand{\bexam}{\begin{example}}
\newcommand{\eexam}{\end{example}}
\newcommand{\bassu}{\begin{assumption}}
\newcommand{\eassu}{\end{assumption}}
\newcommand{\bcoro}{\begin{corollary}}
\newcommand{\ecoro}{\end{corollary}}

\newcommand{\benum}{\begin{enumerate}[label = \emph{(\alph*)}]}
\newcommand{\eenum}{\end{enumerate}}

\newcommand{\mgap}{\vspace{.1in}}

\newcommand{\beq}{\begin{equation}}
\newcommand{\eeq}{\end{equation}}

\DeclareMathOperator{\dom}{dom}

\DeclareMathOperator*{\argmin}{argmin}

\DeclareMathOperator{\prox}{\mbox{prox}}

\newtheorem{theorem}{Theorem}[section] 
\newtheorem{lemma}[theorem]{Lemma}
\newtheorem{corollary}[theorem]{Corollary}
\newtheorem{proposition}[theorem]{Proposition}
\newtheorem{remark}[theorem]{Remark}
\newtheorem{definition}[theorem]{Definition}
\newtheorem{assumption}[theorem]{Assumption}
\newtheorem{example}[theorem]{Example}

\newtheorem{algorithm}{Algorithm}   

\usepackage[T1]{fontenc}
\usepackage[utf8]{inputenc}
\usepackage[portuguese,english]{babel} 

\newcommand{\GG}{{\mathcal{G}}} 
\newcommand{\ES}{{\mathcal{S}}} 
\newcommand{\what}[1]{\widehat #1}
\newcommand{\wtil}[1]{\widetilde #1}

\DeclareMathOperator*{\minimize}{minimize} 
\DeclareMathOperator*{\maximize}{maximize} 

\makeatletter
\def\widebreve#1{\mathop{\vbox{\m@th\ialign{##\crcr\noalign{\kern\p@}%
  \brevefill\crcr\noalign{\kern0.1\p@\nointerlineskip}%
  $\hfil\displaystyle{#1}\hfil$\crcr}}}\limits}

\def\brevefill{$\m@th \setbox\z@\hbox{}%
 \hfill\scalebox{0.7}{\rotatebox[origin=c]{90}{(}} \kern4pt $}
\makeatletter

\begin{document}

\title{A relative-error inexact ADMM splitting algorithm for convex optimization with inertial effects}
\author{
    M. Marques Alves
\thanks{
Departamento de Matem\'atica,
Universidade Federal de Santa Catarina,
Florian\'opolis, Brazil, 88040-900 ({\tt maicon.alves@ufsc.br}).
The work of this author was partially supported by CNPq grant 308036/2021-2.
}
\and
  Marina Geremia
\thanks{
Departamento de Matem\'atica,
Universidade Federal de Santa Catarina,
Florian\'opolis, Brazil, 88040-900. Departamento de Ensino, Pesquisa e Extensão, Instituto Federal de Santa Catarina (IFSC) ({\tt marina.geremia@ifsc.edu.br}).}
}
\date{\emph{Dedicated to the memory of Professor Hedy Attouch}}

\maketitle
\pagestyle{plain}

\begin{abstract}
We propose a new relative-error inexact version of the alternating direction method of multipliers (ADMM) for convex optimization.  
We prove the asymptotic convergence of our main algorithm as well as pointwise and ergodic iteration-complexities for
residuals.
We also justify the effectiveness of the proposed algorithm through some preliminary numerical experiments on 
regression problems.
\\
\\ 
  2000 Mathematics Subject Classification: 90C25, 90C06, 49J52.
 \\
 \\
  Key words: Convex optimization, ADMM, inexact, relative-error, inertial algorithms. 
\end{abstract}

\pagestyle{plain}


\section{Introduction}
Consider the convex problem
\begin{align} \lab{eq:prob}
\displaystyle{\minimize_{x\in \HH} f(x)+g(Lx)},
\end{align}
where $f:\HH\to (-\infty, \infty]$ and $g:\GG \to (-\infty, \infty]$ are lower semicontinuous proper convex functions 
and $L:\HH\to \GG$ is a linear operator ($\HH$ and $\GG$ denote finite-dimensional inner product spaces).
Problem \eqref{eq:prob} appears in different contexts in applied mathematics, including optimization, inverse problems, machine learning, among others.

One of the most popular numerical algorithms for solving \eqref{eq:prob} is the alternating direction method of multipliers (ADMM)~\cite{fortin1983augmented, gabay1983applications, glo.mar-sur.fnm75}, which has now attracted a lot of attention from the numerical optimization community 
(see, e.g., \cite{attouch.augmented2009,valkonen.preconditioned2016,bot.acm2019,boy.par.chu-dis.ftml11,bredies.proximal2017,sun.equivalence2021,eck.ber-dou.mp92,eck.yao-rel.mp18,glo.osh.yin-spl.spi16,hager.coap2020,he.convergence2012,lorenz.nonstationay2019,mon.sva-blo.siam13,teboulle.rate2014,svaiter-par.opt2021}).

In this paper we propose and study a new inexact version of the ADMM allowing relative-error criteria for the solution of the second subproblem (which will appear in the formulation of the proposed algorithm)  
and promoting inertial effects on iterations.

\mgap
\noindent
{\bf Organization of the paper.}
The material is organized as follows. 
In Section \ref{sec:mot} we motivate the definition of our main algorithm, review some related works and discuss the main contributions of this paper.
In Section \ref{sec:main}, we present our main algorithm (Algorithm \ref{alg:main}) and some preliminary results that will be needed to study its convergence and iteration-complexity.
In Section \ref{sec:convergence}, we study the asymptotic behavior of Algorithm \ref{alg:main} under different assumptions under the inertial parameters involved in the formulation of the method. The main results are 
Theorems \ref{th:main} and \ref{th:main02}.
In Section \ref{sec:global} we study the iteration-complexities of Algorithm \ref{alg:main}. The main results in this section are Theorems \ref{the:point} and \ref{the:ergodic}.
Numerical experiments will be presented in Section \ref{sec:num}. 
Appendix \ref{sec:app} contains some auxiliary results.

\section{Motivation, related works and contributions}
  \lab{sec:mot}

\noindent
{\bf Motivation.} We first note that \eqref{eq:prob} is clearly equivalent to the separable problem 
\begin{align}\lab{eq:prob2}
\begin{aligned}
&\minimize \;\; f(x) + g(y),\\
&\mbox{subject to}\;\; Lx - y =0.
\end{aligned}
\end{align}

\mgap

\noindent
An iteration of the standard ADMM~\cite{eck.ber-dou.mp92} for solving \eqref{eq:prob2} can be described as follows: given a starting point
$(y_0,z_0)\in \GG^2$ and a regularization parameter $\gamma>0$, iterate for $k\geq 0$:
\begin{align}
\lab{eq:sadmm01}
&x_{k+1} \in \argmin_{x\in \HH}\,\left\{ f(x) + \inner{z_k}{Lx - y_k} + \dfrac{\gamma}{2}\normq{Lx - y_k}\right\},\\[2mm]
\lab{eq:sadmm02}
&y_{k+1} \in \argmin_{y\in \GG}\,\left\{g(y) + \inner{z_k}{Lx_{k+1} -y }
   +\dfrac{\gamma}{2}\normq{Lx_{k+1}-y}\right\},\\[2mm]
\lab{eq:sadmm03}
&z_{k+1} = z_k + \gamma\left(Lx_{k+1} - y_{k+1}\right).
\end{align}

\mgap

\noindent
We consider here the case in which \eqref{eq:sadmm01} can be solved exactly and, on the other hand, \eqref{eq:sadmm02} is supposed to be solved only approximately by some other (inner) algorithm, like, for instance, CG or BFGS, depending on the particular structure of the function $g(\cdot)$ in \eqref{eq:prob2}.

With this in mind, we will introduce a notion of relative-error approximate solution for \eqref{eq:sadmm02} (more details will be given on Section \ref{sec:main}). 
To this end, first note that \eqref{eq:sadmm02} is an instance of the general family of minimization problems
\begin{align}\lab{eq:sadmm04}
\minimize_{y\in \GG}\left\{g(y) + \inner{z}{Lx - y }
   +\dfrac{\gamma}{2}\normq{Lx - y}\right\},
\end{align}
where $x\in \HH$,  $z\in \GG$ and $\gamma>0$ are given (in the case of \eqref{eq:sadmm02}, we have
$(x, z) = (x_{k+1}, z_k)$).
Moreover, since the function $g(\cdot)$ is convex, we have that \eqref{eq:sadmm04} is also equivalent to the
inclusion/equation system for the pair $(y,v)$:
\begin{align}\lab{eq:ie_sys}
\begin{aligned}
\begin{cases}
v\in \partial g(y),\\[2mm]
v - z + \gamma(y - Lx) = 0.
\end{cases}
\end{aligned}
\end{align}

A formal definition of approximate (inexact) solution of \eqref{eq:ie_sys} (or, equivalently, \eqref{eq:sadmm04}) will be given in Definition \ref{def:app_sol} in Section \ref{sec:main}; such a notion of approximate solution will allow for errors in both the inclusion and the equation in \eqref{eq:ie_sys}.

\mgap
\noindent
{\bf The extended-solution set.} The Fenchel dual of \eqref{eq:prob} is
\begin{align} \lab{eq:d_prob}
\displaystyle{\maximize_{z\in \GG} -f^*(-L^*z) - g^*(z)},
\end{align}
where $f^*:\HH\to (-\infty,\infty]$ and $g^*:\GG\to (-\infty,\infty]$ denote the Fenchel conjugates of $f$ and g, respectively, and $L^*:\GG\to \HH$ denotes the adjoint operator of $L$.
Under standard regularity conditions~\cite{bau.com-book} on $f,g$ and $L$ it is well-known that \eqref{eq:prob}
and \eqref{eq:d_prob} are, respectively, equivalent to the (monotone) inclusions 
\begin{align} \lab{eq:mono}
0 \in \partial f(x) + L^*\partial g(Lx),
\end{align}
and
\begin{align} \lab{eq:d_mono}
0 \in -L\partial f^*(-L^*z) + \partial g^*(z).
\end{align}

We make the blanket assumption:
\bassu
 \lab{ass:blanket}
For the function $f$ and the operator $L$ as in \eqref{eq:prob}, the following holds:
\[
 \partial (f^*\circ -L^*) = -L\circ\partial f^*\circ -L^*.
\]
\eassu
Several sufficient conditions for Assumption \ref{ass:blanket} to hold true can be found, e.g., in \cite{bau.com-book}.
We will also consider an extended-solution set $\mathcal{S}$, 
attached to the pair of inclusions \eqref{eq:mono}--\eqref{eq:d_mono}, defined as
\begin{align} \lab{eq:exts}
\mathcal{S} =\set{(z,w)\in \GG^2\;\;|\;\; -w\in \partial(f^*\circ -L^*)(z)\;\;\mbox{and}\;\; w\in \partial g^*(z)}.
\end{align}

Under Assumption \ref{ass:blanket}, it is easy to check that if $(z,w)\in \mathcal{S}$, then it follows that
there exists $x\in \HH$ such that $x\in \partial f^*(-L^*z)$, $w = Lx$ and 
$x$ and $z$ are solutions of \eqref{eq:mono} and \eqref{eq:d_mono}, respectively.

\mgap

Throughout this work we will assume the following.
\bassu
 \lab{ass:nonempty}
We assume the extended solution set $\ES$ as in \eqref{eq:exts} is nonempty.
\eassu

\noindent
{\bf Inertial algorithms.} Iterative algorithms with inertial effects for monotone inclusions (and related topics in 
optimization, saddle-point, equilibrium problems, etc) were first proposed in the seminal paper~\cite{alv.att-iner.svva01}
and subsequently developed in various directions of research by different authors and research groups (see, e.g., 
\cite{alv.eck-admm.coap20, att-fast.minimax2021,att-convergence.mp2020,att-inertial.jota2018,attouch.peypouquet-convergence.mp2019,bot.acm2019,combettes-quasi.sopt2017} and references therein).
Basically, the main idea consists in at a current iterate, say $p_k$, produce an ``inertial effect'' by a simple extrapolation:
\[
 \what p_k = p_k + \alpha_k(p_k - p_{k-1}),
\]
where $\alpha_k\geq 0$, and then generate the next iterate $p_{k+1}$ from $\what p_k$ instead of $p_k$ 
(see \eqref{eq:extrap}--\eqref{eq:extrap3} below). 
Our main algorithm, namely Algorithm \ref{alg:main}, will benefit from inertial effects on the iteration; see the comments and remarks following Algorithm \ref{alg:main} for more discussions regarding the effects of inertia. 

\mgap

\noindent
{\bf Main contributions.} 
We present a theoretical (asymptotic and iteration-complexity analysis) and computational study of a partially inexact  
ADMM splitting algorithm for solving \eqref{eq:prob}.
Our main algorithm, namely Algorithm \ref{alg:main} below, benefits from the addition of inertial effects; 
see \eqref{eq:extrap} and \eqref{eq:extrap3}.
The convergence analysis is presented in Theorem \ref{th:main}, to which the proof incorporates some elements 
of \cite{alv.eck-admm.coap20} and \cite{svaiter-par.opt2021}.
We also obtained iteration-complexities for the proposed algorithm by showing pointwise $O(1/\sqrt{k})$ and ergodic $O(1/k)$ global convergence rates for residuals; see Theorems \ref{the:point} and \ref{the:ergodic} below.
We justify the effectiveness of our main algorithm through the realization of numerical experiments on 
the LASSO problem (see Section \ref{sec:num}). 

\mgap
\noindent
{\bf Related works.} A partially inexact ADMM splitting algorithm was recently proposed and studied 
in \cite{svaiter-par.opt2021}.
Paper~\cite{melo2019partially} proposes a partially inexact ADMM for which the first subproblem is supposed to be solved inexactly. The analysis of the main algorithm in~\cite{melo2019partially} is performed by viewing it as a special instance of a non-Euclidean version of the hybrid proximal extragradient method~\cite{mon.sva-hpe.siam10}. 
In contrast to this, analogously to~\cite{svaiter-par.opt2021}, our main algorithm (Algorithm \ref{alg:main} below)
assumes the second subproblem is solved inexactly.
Moreover, since \eqref{eq:wood}--\eqref{eq:wood2} below also allows for errors in $\partial g$, the error criterion we propose here is potentially more flexible than the corresponding one in~\cite{melo2019partially}.
Other relative-error inexact versions of ADMM were also previously studied in~\cite{xie2018inexact,xie2017inexact},
but we notice that the convergence results were restricted to the analysis of the dual sequences.  
We also mention that the relative-error inexact variants of the ADMM from~\cite{alv.eck-admm.coap20,eck.yao-rel.mp18} only apply to \eqref{eq:prob} in the particular case of $L=I$, and, additionally, these variants assume the first subproblem to be solved inexactly with the error condition verified only a-posteriori, that is, only after the computation of second subproblem's solution.

\mgap

\noindent
{\bf General notation.} We denote by $\HH$ and $\GG$ finite-dimensional inner product spaces with inner product 
and induced norm denoted, respectively, by $\inner{}{}$ and $\norm{\cdot}=\sqrt{\inner{\cdot}{\cdot}}$.
For any set $\mathcal{X}$ we denote by $X^n$ the $n$-product $\mathcal{X}\times \cdots \times \mathcal{X}$.
In $\GG^2$, we will consider the inner product and induced norm defined, respectively, by 
\begin{align} \lab{def:inner_g3}
\inner{p}{p'}_\gamma := \dfrac{1}{\gamma} \inner{z}{z'} + \gamma \inner{w}{w'}\;\;\mbox{and}\;\;
\normq{p}_\gamma := \inner{p}{p}_\gamma,
\end{align}
where $p=(z,w), p'=(z',w')\in \GG^2$ and $\gamma>0$. More precisely, for $p=(z, w)\in \GG^2$, the norm of 
$p$ is
\begin{align} \lab{def:norm_g3}
\normq{p}_\gamma = \dfrac{1}{\gamma}\normq{z} + \gamma \normq{w}.
\end{align}
An extended-real valued function $f:\HH\to (-\infty, \infty]$ is said to be \emph{convex} whenever
$f(\lambda x + (1-\lambda)y)\leq \lambda f(x)+(1-\lambda)f(y)$ for all $x,y\in \HH$ and $\lambda\in (0,1)$,
and $f$ is \emph{proper} if its \emph{effective domain}, denoted by $\dom f$, is nonempty.
The \emph{Fenchel conjugate} of a proper function $f:\HH\to (-\infty,\infty]$ is $f^*:\HH\to (-\infty,\infty]$, defined
at any $u\in \HH$ by $f(u)=\sup_{x\in \HH}\,\{\inner{x}{u}-f(x)\}$.
The \emph{$\varepsilon$-subdifferential} and the \emph{subdifferential} of a convex function 
$g:\HH\to (-\infty,\infty]$ at $x\in \HH$ are defined as
$\partial_\varepsilon g(x) := \{u\in \HH\;|\; g(y)\geq g(x)+\inner{u}{y-x}-\varepsilon\quad \forall y\in \HH\}$ and 
$\partial g(x) := \partial_0 g(x)$, respectively. For additional details on standard notations and definitions of convex analysis we refer the reader to the references~\cite{bau.com-book,rock-ca.book}.

\section{The main algorithm and some preliminary results}
 \lab{sec:main}
 
Consider the minimization problem \eqref{eq:prob}, i.e.,

\begin{align} \lab{eq:probx}
\displaystyle{\minimize_{x\in \HH} f(x)+g(Lx)},
\end{align}
where $f:\HH\to (-\infty,\infty]$ and $g:\GG \to (-\infty,\infty]$ are lower semicontinuous proper convex functions 
and $L:\HH\to \GG$ is a linear operator between finite-dimensional inner product spaces $\HH$ and $\GG$.

In this section we present our main algorithm, namely Algorithm \ref{alg:main} below. This is a partially inexact (the second block is allowed to be solved inexactly) ADMM with relative-error criterion for the second subproblem. 
Recall the extended solution set $\mathcal{S}$ as in \eqref{eq:exts} and Assumptions \ref{ass:blanket} and \ref{ass:nonempty}.
The three technical lemmas \ref{lm:signal}, \ref{lem:psigma}, \ref{lm:ineq.delta} and \ref{lem:dina} will be used in the subsequent section. 

Before presenting our main algorithm, as we discussed in the Introduction, we have to formalize the notion of inexact 
solution that will be used to compute approximate solution for the second subproblem. Recall that the second subproblem of the 
standard ADMM (see \eqref{eq:sadmm02}) belongs to the general family of minimization problems \eqref{eq:sadmm04}, which is, in particular, equivalent to the inclusion/equation system \eqref{eq:ie_sys} for the pair
$(y,v)$, i.e.,

\begin{align}\lab{eq:ie_sys2}
\begin{aligned}
\begin{cases}
v\in \partial g(y),\\[2mm]
v - z + \gamma(y - Lx) = 0.
\end{cases}
\end{aligned}
\end{align}

\bdefi[$\sigma$-approximate solution of \eqref{eq:sadmm04}]
 \lab{def:app_sol}
For $x\in \HH$, $(\widehat z, \widehat y)\in \GG^2$ and $\gamma > 0$, a triple $(\widetilde y, v, \varepsilon)\in \GG\times \GG\times \R_+$ is said to be a $\sigma$-approximate solution of \eqref{eq:sadmm04} \emph{(}or, equivalently, of \eqref{eq:ie_sys}\emph{)} 
at $(x, \widehat z, \widehat y)$ 
if $\sigma\in [0,1)$ and
\begin{align}\lab{eq:app_sol}
\begin{aligned}
\begin{cases}
v\in \partial_\varepsilon g(\widetilde y),\\[2mm]
v - \widehat z + \gamma(\widetilde y - Lx) =: e\\[2mm]
\normq{e} + 2\gamma\varepsilon \leq \sigma^2 \min\left\{\gamma^2\normq{Lx-\widehat y}, \normq{v - \widehat z}\right\}.
\end{cases}
\end{aligned}
\end{align}
We will also write 
\begin{align*}
\widetilde y\overset{\sigma}{\approx} \argmin_{y\in \GG}\left\{g(y) + \inner{\widehat z}{Lx - y }
   +\dfrac{\gamma}{2}\normq{Lx - y}\right\}
\end{align*}
meaning that there exists $(v,\varepsilon)$ such that $(\widetilde y, v, \varepsilon)$ satisfies \eqref{eq:app_sol}.
\edefi

\mgap

\noindent
We now make some remarks regarding Definition \ref{def:app_sol}:
\bitem
\item[(i)] Note that if $\sigma=0$ in \eqref{eq:app_sol}, then it follows that $e=0$ and $\varepsilon=0$, which is to say that
the pair $(\widetilde y, v)$ satisfies the inclusion/equation system \eqref{eq:ie_sys} (recall that $\partial_0 g= \partial g$) and, in particular, $\widetilde y$ is an exact solution of \eqref{eq:sadmm04}.
\item[(ii)] The error criterion for \eqref{eq:ie_sys} as in \eqref{eq:app_sol} belongs to the class of \emph{relative-error} criteria for proximal-type algorithms. Different variants of such error conditions have been employed for computing approximate solution for (sub) problems for a wide range of algorithms in monotone inclusions, convex optimization, 
saddle-point problems, etc 
(see, e.g., \cite{alv.eck-admm.coap20, eck.yao-rel.mp18, mon.sva-hpe.siam10,sol.sva-hpe.svva99,sol.sva-hyb.jca99,sol.sva-uni.nfao01}).
\item[(iii)] The error criterion \eqref{eq:app_sol} will be used to compute approximate solutions in step 3 of our main algorithm, namely Algorithm \ref{alg:main} below (see \eqref{eq:wood}--\eqref{eq:wood2}).
\item[(iv)] As an illustrative example, consider the special case of the LASSO problem~\cite{tib-las.jrss96}
\begin{align}\lab{eq:lasso}
\min_{x\in \R^d}\,\left\{\dfrac{1}{2}\normq{Ax - b}+\nu\norm{x}_1\right\},
\end{align}
where $A\in \R^{n\times d}$, $b\in \R^n$ and $\nu>0$. Problem \eqref{eq:lasso} is clearly a special instance of \eqref{eq:prob} in which $L:=I$, $f(x):=\nu\norm{x}_1$ and $g(x):=(1/2)\normq{Ax-b}$ (see also Section \ref{sec:num} below).
In this case, our inclusion/equation system \eqref{eq:ie_sys} clearly reduces to
\begin{align*}
v = A^*(Ay -b),\qquad v - z + \gamma(y - x) = 0,
\end{align*}
or, in other words, in this special case, \eqref{eq:ie_sys} is equivalent to the linear system (operator equation)
\begin{align*}
\Big(A^*A + \gamma I\Big) y = A^*b + z + \gamma x.
\end{align*}
The latter linear system can be solved by the CG algorithm~\cite{nocedal.book}, where $e$ as in \eqref{eq:app_sol} will simply denote the residual of the system and the inequality in \eqref{eq:app_sol} can be used as a stopping criterion for 
CG.
\eitem

\mgap

Next is our main algorithm.

\mgap

%
%
\noindent
\fbox{
\begin{minipage}[h]{6.6 in}
\begin{algorithm} \label{alg:main}
{\bf An inexact inertial ADMM algorithm for solving \eqref{eq:prob}}
\end{algorithm}
\begin{itemize}
\item[(0)] Let $(z_0, y_0) = (z_{-1}, y_{-1})\in \GG^2$ and let 
$\alpha,\sigma \in [0, 1)$, $\tau\in (0, 1)$ and $\gamma>0$. Set $k=0$.
\item [(1)] Choose $\alpha_k\in [0,\alpha]$ and let
\begin{align}
%
\lab{eq:extrap}
& \what z_k = z_k + \alpha_k(z_k - z_{k-1}),\\
\lab{eq:extrap3}
& \what y_k = y_k + \alpha_k(y_k - y_{k-1}).
%
\end{align}
\item [(2)] Compute 
\begin{align}\lab{eq:borwein}
\hspace{-3.2cm}x_k \in \argmin_{x\in \HH}\,\left\{f(x)+\inner{\what z_k}{Lx - \what y_k} + \dfrac{\gamma}{2}\normq{Lx - \what y_k}\right\}.
\end{align}
 \item[(3)] Compute
\begin{align}\lab{eq:light}
\widetilde y_k \overset{\sigma}{\approx} \argmin_{y\in \GG}\,
   \left\{g(y)+\inner{\what z_k}{Lx_k - y} + \dfrac{\gamma}{2}\normq{Lx_k - y}\right\}
\end{align}
at $(x_k, \widehat z_k, \widehat y_k)$
 in the sense of Definition \ref{def:app_sol}, i.e., compute $(\widetilde y_k, v_k, \varepsilon_k)\in \GG\times \GG\times \R_+$ such that
 \begin{align}
 \lab{eq:wood}
 & v_k \in \partial_{\varepsilon_k} g(\widetilde y_k),\\[2mm]
 \lab{eq:woodx}
 & v_k - \widehat z_k + \gamma(\widetilde y_k - L x_k) =: e_k,\\[2mm]
 \lab{eq:wood2}
 & \normq{e_k}
  +2\gamma\varepsilon_k
  \leq  \sigma^2 \min\left\{\gamma^2\normq{Lx_k - \what y_k}, \normq{v_k - \widehat z_k}\right\}.
 \end{align}
\item[(4)] Set
\begin{align} 
\lab{eq:belchior}
& z_{k+1} = \what z_k + \tau \gamma(Lx_k - \widetilde y_k),\\[2mm]
\lab{eq:belchior2}
& y_{k+1} = (1 - \tau)\what y_k + 
    \frac{\tau}{\gamma} (\widehat z_k + \gamma L x_k - v_k),
 \end{align}
 $k = k+1$ and go to step 1.
\end{itemize}
\end{minipage}
} 
%

\mgap
\mgap

\noindent
We now make some remarks concerning Algorithm \ref{alg:main}:
\bitem
\item[(i)] Algorithm \ref{alg:main} is specially designed for instances of \eqref{eq:prob} in which \eqref{eq:borwein} has a closed-form solution, i.e., for problems in which \eqref{eq:borwein} is easy to solve.
In this regard, one example of interest is when $f(\cdot) = \norm{\cdot}_1$ and $L=I$, in which case \eqref{eq:borwein} has a unique solution given explicitly by $x_k = \prox_{\gamma^{-1}\|\cdot\|_1}(\what y_k - \gamma^{-1}\what z_k)$.
On the other hand, we assume that the computation of $\wtil y_k$ as in \eqref{eq:light} demands the use of an (inner) algorithm, which the choice of depends on the particular structure of the function $g$, and, in this case, one can use 
\eqref{eq:wood}--\eqref{eq:wood2} as a stopping criterion for the inner algorithm of choice.
\item[(ii)] Recall that we discussed in the Introduction (see ``Related works'') other ADMM-type algorithms related to Algorithm \ref{alg:main}.
\item[(iii)] The main results on the convergence and iteration-complexity of Algorithm \ref{alg:main} are Theorems \ref{th:main}, \ref{th:main02}, \ref{the:point} and \ref{the:ergodic} below.  Numerical experiments will be presented and discussed in Section \ref{sec:num}.
\item[(iv)] The role of the parameter $0< \tau<1$ is to introduce (under) relaxation in the iterative process; see
\eqref{eq:belchior} and \eqref{eq:belchior2}.
\item[(v)] We will also need the sequences $(\breve{z}_k)$ and $(\breve{y}_k)$, where, for all $k\geq 0$,
\begin{align} \lab{eq:latam}
\breve{z}_k := \what z_k + \gamma(Lx_k - \widetilde y_k),\quad
\breve{y}_k := \frac{1}{\gamma} (\widehat z_k + \gamma L x_k - v_k).
 \end{align}
Note that $\breve{z}_k = z_{k+1}$ and $\breve{y}_k = y_{k+1}$ if we set $\tau = 1$ in \eqref{eq:belchior} and \eqref{eq:belchior2},
respectively.
\eitem
\mgap
Next we present four technical lemmas -- Lemmas \ref{lm:signal}, \ref{lem:psigma}, \ref{lm:ineq.delta} and 
\ref{lem:dina} --, which will be useful in the subsequent sections. 

\mgap

\blemm \lab{lm:signal}
Consider the sequences evolved by \emph{Algorithm \ref{alg:main}}, let $\ES$ be as in \eqref{eq:exts} and let $(\breve{z}_k)$
and $(\breve{y}_k)$ be as in \eqref{eq:latam}. Define
\begin{align}\lab{eq:def.pk}
 p_k = (z_k, y_k), \quad \what p_k=(\what z_k, \what y_k), \quad \breve{p}_k = (\breve{z}_k, \breve{y}_k)\;\;\emph{and}\;\;
 \wtil p_k = (v_k, L x_k) \qquad \forall k\geq 0.
\end{align}
\begin{enumerate}[label = \emph{(\alph*)}]
\item \lab{ite:signalm}
For all $k\geq 0$,
\begin{align*}
p_{k+1} = (1 - \tau)\what p_k + \tau \breve{p}_k.
\end{align*}
\item \lab{ite:signaln}
For all $k\geq 0$,
\begin{align*}
 -Lx_k \in \partial \left(f^*\circ -L^*\right)(z'_k),
\end{align*}
where
\begin{align} \lab{def:zkk}
z'_k := \what z_k + \gamma\left(L x_k - \what y_k\right). 
\end{align}
\item \lab{ite:signalo}
For all $k\geq 0$,
\begin{align*}
 \breve{p}_k - \widehat p_k = \left(\gamma(L x_k - \widetilde y_k), \frac{1}{\gamma}\left(z'_k - v_k\right)\right).
\end{align*}
\item \lab{ite:signalp} 
For all $k\geq 0$ and $p=(z, w)\in \ES$,
\begin{align*}
\inner{\breve{p}_k - \what p_k}{p - \wtil p_k}_\gamma\geq -\varepsilon_k.
\end{align*}
\item \lab{ite:signalq} 
For all $k\geq 0$,
\begin{align*}
\what p_k = p_k +\alpha_k(p_k - p_{k-1}).
\end{align*}
\end{enumerate}
\elemm
\mgap
\bproo
\iteref{ite:signalm} This result is a direct consequence of the definitions of $p_k$, $\what p_k$ and $\breve{p}_k$
as in \eqref{eq:def.pk} -- see also \eqref{eq:latam} -- combined with \eqref{eq:belchior} and \eqref{eq:belchior2}.

\mgap

\iteref{ite:signaln} First note that from \eqref{eq:borwein} and \eqref{def:zkk}, we obtain
$
0\in \partial f(x_k) + L^* z'_k,
$
or, equivalently, $-L^* z'_k \in \partial f(x_k)$. As $\left(\partial f\right)^{-1}=\partial f^*$, the latter inclusion is also 
equivalent to 
$x_k \in \partial f^*(-L^* z'_k)$, which in turn yields 
$-Lx_k \in -L\partial f^*(-L^* z'_k)$, which by Assumption \ref{ass:blanket} gives item \iteref{ite:signaln}.

\mgap

\iteref{ite:signalo} This follows easily from \eqref{eq:latam} -- \eqref{def:zkk} and some simple algebraic manipulations.

\mgap

\iteref{ite:signalp} As $\left(\partial_{\varepsilon_k} g\right)^{-1} = \partial_{\varepsilon_k} g^*$ -- see, e.g., \cite[p. 85, Theorem 2.4.4(iv)]{Zal02} --, we have
that \eqref{eq:wood} is equivalent to the inclusion
\begin{align}\lab{eq:s_inc}
\widetilde y_k\in \partial_{\varepsilon_k} g^*(v_k).
\end{align}

As $p = (z, w)\in \ES$, according to the definition of $\ES$ in \eqref{eq:exts}, 
we have $-w \in \partial \left(f^*\circ -L^*\right)(z)$ 
and $w\in \partial g^*(z)$. The latter inclusions combined with item \iteref{ite:signaln} and \eqref{eq:s_inc} 
and the monotonicity of $\partial(f^*\circ -L^*)$ and $\partial g^*$ 
yield
\begin{align}\lab{eq:sunday}
\inner{z'_k - z}{w - Lx_k}\geq 0\;\;\mbox{and}\;\;\inner{z - v_k}{w - \widetilde y_k}\geq -\varepsilon_k.
\end{align}
Now using \eqref{eq:def.pk}, item \iteref{ite:signalo}, \eqref{eq:sunday}  and the definition of $\inner{\cdot}{\cdot}_\gamma$ as in \eqref{def:inner_g3} we find
\begin{align*}
\inner{\breve{p}_k - \what p_k}{p - \wtil p_k}_\gamma & = 
 \dfrac{1}{\gamma}\Inner{\gamma(L x_k - \widetilde y_k)}{z - v_k} 
 + \gamma \Inner{\frac{1}{\gamma}\left(z'_k - v_k\right)}{w - L x_k}\\[2mm]
& = \inner{L x_k - \widetilde y_k}{z - v_k} + \inner{z'_k - z}{w - Lx_k} + \inner{z - v_k}{w - Lx_k}\\[2mm]
& = \inner{z - v_k}{w - \tilde y_k} + \inner{z'_k - z}{w - Lx_k}\\[2mm]
&\geq -\varepsilon_k.
\end{align*}

\mgap

\iteref{ite:signalq} This follows directly from \eqref{eq:extrap}, \eqref{eq:extrap3} and the definition of $\what p_k$ as in \eqref{eq:def.pk}.
\eproo

\mgap

\blemm \lab{lem:psigma}
Consider the sequences evolved by \emph{Algorithm \ref{alg:main}} and
let $(\breve{p}_k)$, $(\what p_k)$ and $(\wtil p_k)$ be as in \eqref{eq:def.pk}. 
For all $k\geq 0$,
\begin{enumerate}[label = \emph{(\alph*)}]
\item \lab{ite:psigmam}
$\normq{\wtil p_k - \breve{p}_k}_\gamma = \frac{1}{\gamma}\Big(\normq{e_k} + 
      \normq{v_k - \what z_k}\Big).$
\item \lab{ite:psigman}
$\norm{\breve{p}_k - \what p_k}_\gamma \leq 2\norm{\wtil p_k - \what p_k}_\gamma.$
\end{enumerate}
\elemm
\bproo
\iteref{ite:psigmam} Direct use of \eqref{def:norm_g3}, \eqref{eq:woodx}, \eqref{eq:latam} and \eqref{eq:def.pk} gives
\begin{align*}
\nonumber
\normq{\wtil p_k - \breve{p}_k}_\gamma &= \frac{1}{\gamma}\normq{v_k - \breve{z}_k}
    + \gamma \normq{L x_k - \breve{y}_k}\\[2mm]
 \nonumber
 & = \frac{1}{\gamma}\normq{\underbrace{v_k - \what z_k + \gamma(\widetilde y_k - Lx_k)}_{e_k}} 
   + \gamma \normq{Lx_k - \left[\gamma^{-1} (\widehat z_k + \gamma L x_k - v_k)\right]}\\[2mm]
  & = \frac{1}{\gamma}\Big(\normq{e_k} + 
      \normq{v_k - \what z_k}\Big).
\end{align*}

\mgap

\iteref{ite:psigman} In view of \eqref{eq:wood2}, \eqref{eq:def.pk} and item (a),
\begin{align*}
\normq{\wtil p_k - \breve{p}_k}_\gamma & = \frac{1}{\gamma}\left(\normq{e_k} + 
      \normq{v_k - \what z_k}\right)\\[2mm]
& \leq \dfrac{1}{\gamma}\left(\sigma^2 \normq{\gamma(Lx_k - \what y_k)} + \normq{v_k - \what z_k}\right)\\[2mm]
&\leq \dfrac{1}{\gamma}\normq{v_k - \what z_k} + \gamma\normq{Lx_k - \what y_k}\\[2mm]
& = \normq{\wtil p_k - \what p_k}_\gamma.
\end{align*}
Hence, using the triangle inequality,
\begin{align*}
\norm{\breve{p}_k - \what p_k}_\gamma \leq \norm{\wtil p_k - \breve{p}_k}_\gamma + \norm{\wtil p_k - \what p_k}_\gamma \leq 2\norm{\wtil p_k - \what p_k}_\gamma.
\end{align*}
\eproo

\mgap

\blemm \lab{lm:ineq.delta}
Consider the sequences evolved by \emph{Algorithm \ref{alg:main}} and let $(\what p_k)$, $(\breve{p}_k)$
and $(\wtil p_k)$ be as in \eqref{eq:def.pk}.
\benum
\item \lab{lm:ineq.deltam}
For all $k\geq 0$ and $p \in \GG^2$,
\[
\normq{p-\what p_k}_\gamma - \normq{p - \breve{p}_k}_\gamma = 
 \normq{\wtil p_k - \what p_k}_\gamma - \normq{\wtil p_k - \breve{p}_k}_\gamma
 + 2\inner{\breve{p}_k - \what p_k}{p - \wtil p_k}_\gamma.
\]
\item \lab{lm:ineq.deltan}
For all $k\geq 0$ and $p \in \GG^2$,
\[
\normq{p-\what p_k}_\gamma - \normq{p - \breve{p}_k}_\gamma\geq \gamma(1 - \sigma^2) \normq{L x_k - \what y_k} 
+ 2\left[\varepsilon_k + \inner{\breve{p}_k - \what p_k}{p - \wtil p_k}_\gamma\right].
\]
\item \lab{lm:ineq.deltao}
For all $k\geq 0$ and $p = (z, w) \in \ES$,
\[
\normq{p-\what p_k}_\gamma - \normq{p - \breve{p}_k}_\gamma\geq \gamma(1 - \sigma^2) \normq{L x_k - \what y_k}.
\]
\item \lab{lm:ineq.deltap}
For all $k\geq 0$ and $p = (z, w) \in \ES$,
\[
\normq{p-\what p_k}_\gamma - \normq{p - p_{k+1}}_\gamma\geq 
\tau\gamma(1-\sigma^2) \normq{Lx_k - \what y_k} + (1-\tau)\tau\normq{\breve{p}_k - \what p_k}_\gamma.
\]
\item \lab{lm:ineq.deltaq}
For all $k\geq 0$ and $p = (z, w) \in \ES$,
\begin{align*}
\normq{p-\what p_k}_\gamma - \normq{p - p_{k+1}}_\gamma&\geq 
\tau(1 - \tau)(1 - \sigma)^2\normq{\wtil p_k - \what p_k}_\gamma\\[2mm]
& \geq \dfrac{(1 - \tau)(1 - \sigma)^2}{4\tau}\normq{p_{k+1} - \what p_k}_\gamma.
\end{align*}
\eenum
\elemm
\bproo
\iteref{lm:ineq.deltam} The desired result follows directly from the well-known identity 
$\normq{a-b}_\gamma - \normq{a-c}_\gamma=\normq{d-b}_\gamma - \normq{d-c}_\gamma+2\inner{c-b}{a-d}_\gamma$
with $a=p$, $b=\what p_k$, $c = \breve{p}_k$ and $d=\wtil p_k$.

\mgap

\iteref{lm:ineq.deltan} Using \eqref{def:norm_g3} and the definitions of $(\wtil p_k)$ and $(\what p_k)$ as in \eqref{eq:def.pk} we get
\begin{align*}
\normq{\wtil p_k - \what p_k}_\gamma  = \frac{1}{\gamma}\normq{v_k - \what z_k} + \gamma \normq{Lx_k - \what y_k},
\end{align*}
which in turn combined with Lemma \ref{lem:psigma}\iteref{ite:psigmam} yields
\begin{align*}
\normq{\wtil p_k - \what p_k}_\gamma - \normq{\wtil p_k - \breve{p}_k}_\gamma = \gamma \normq{Lx_k - \what y_k}
   -\dfrac{1}{\gamma} \normq{e_k}.
\end{align*}
From \eqref{eq:wood2}, item (a), the latter identity and some algebraic manipulations,
\begin{align*}
\normq{p-\what p_k}_\gamma - \normq{p - \breve{p}_k}_\gamma &= 
\gamma \normq{Lx_k - \what y_k} - \dfrac{1}{\gamma} \normq{e_k}  + 2\inner{\breve{p}_k - 
  \what p_k}{p - \wtil p_k}_\gamma\\[2mm]
& = \gamma \normq{Lx_k - \what y_k} - \dfrac{1}{\gamma}\left(\normq{e_k} + 2\gamma\varepsilon_k\right) 
 + 2\left[\varepsilon_k + \inner{\breve{p}_k - 
  \what p_k}{p - \wtil p_k}_\gamma\right]\\[2mm]
&\geq \gamma(1 - \sigma^2)\normq{Lx_k - \what y_k} + 2\left[\varepsilon_k + \inner{\breve{p}_k - 
  \what p_k}{p - \wtil p_k}_\gamma\right],
\end{align*}
which finishes the proof of \iteref{ite:psigman}.

\mgap

\iteref{lm:ineq.deltao} This is a direct consequence of Lemma \ref{lm:signal}\iteref{ite:signalp} and item \iteref{lm:ineq.deltan} above. 

\mgap

\iteref{lm:ineq.deltap} Using Lemma \ref{lm:signal}\iteref{ite:signalm} and the identity 
$\normq{(1-\tau) a + \tau b}_\gamma = (1-\tau)\normq{a}_\gamma + \tau\normq{b}_\gamma -(1-\tau)\tau\normq{a - b}_\gamma$
with $a = p - \what p_k$ and $b = p - \breve{p}_k$, we obtain
\begin{align*}
\nonumber
\normq{p - p_{k+1}}_\gamma &= \normq{(1-\tau)(p - \what p_k) + \tau(p - \breve{p}_k)}_\gamma\\[2mm]
  & = (1-\tau)\normq{p - \what p_k}_\gamma + \tau\normq{p - \breve{p}_k}_\gamma - (1-\tau)\tau\normq{\breve{p}_k - \what p_k}_\gamma.
\end{align*}
Now by multiplying the inequality in item \iteref{lm:ineq.deltao} by $\tau>0$, using the latter identity and some simple algebraic manipulations, we find
the desired result.

\mgap

\iteref{lm:ineq.deltaq} Note first that using \eqref{eq:wood2} and the triangle inequality, we find
\begin{align*}
\norm{v_k - \what z_k}&\leq \norm{\underbrace{v_k - \what z_k + \gamma(\wtil y_k - Lx_k)}_{e_k}} + \norm{\gamma(\wtil y_k - Lx_k)}\\
   &\leq \sigma \norm{v_k - \what z_k} +  \norm{\gamma(\wtil y_k - Lx_k)},
\end{align*}
so that
\begin{align*}
\norm{\gamma(\wtil y_k - Lx_k)} \geq (1 - \sigma)\norm{v_k - \what z_k},
\end{align*}
which in turn combined with \eqref{def:norm_g3},  \eqref{eq:latam} and \eqref{eq:def.pk} yields
\begin{align} \lab{eq:attc02}
\nonumber
\normq{\breve{p}_k - \what p_k}_\gamma & = \dfrac{1}{\gamma}\normq{\breve{z}_k - \what z_k} + \gamma \normq{\breve{y}_k - \what y_k}\\[2mm]
\nonumber
 &\geq \dfrac{1}{\gamma}\normq{\breve{z}_k - \what z_k}\\[2mm]
\nonumber
& = \dfrac{1}{\gamma}\normq{\gamma(\wtil y_k - Lx_k)}\\[2mm]
& \geq \dfrac{1}{\gamma}(1 - \sigma)^2\normq{v_k - \what z_k}.
\end{align}

Now using \eqref{def:inner_g3}, \eqref{eq:def.pk}, item \iteref{lm:ineq.deltap} above and \eqref{eq:attc02}, 
\begin{align} \lab{eq:ivan}
\nonumber
 \normq{p-\what p_k}_\gamma - \normq{p - p_{k+1}}_\gamma&\geq 
\tau\gamma(1-\sigma^2) \normq{Lx_k - \what y_k} + (1-\tau)\tau\normq{\breve{p}_k - \what p_k}_\gamma\\[2mm]
\nonumber
  &\geq  \tau\gamma(1-\sigma^2) \normq{Lx_k - \what y_k} + (1-\tau)\tau\dfrac{1}{\gamma}(1 - \sigma)^2\normq{v_k - \what z_k}
\\[2mm]
\nonumber
& \geq \tau(1 - \tau)(1 - \sigma)^2
\Big[\dfrac{1}{\gamma}\normq{v_k - \what z_k} + \gamma\normq{Lx_k - \what y_k}\Big]\\[2mm]
& = \tau(1 - \tau)(1 - \sigma)^2\normq{\wtil p_k - \what p_k}_\gamma.
\end{align}
To prove the second inequality, one can use Lemmas \ref{lm:signal}\iteref{ite:signalm} and \ref{lem:psigma}\iteref{ite:psigman} to conclude that $\norm{p_{k+1} - \what p_k}_\gamma \leq 2\tau \norm{\wtil p_k - \what p_k}_\gamma$ and then apply it in
\eqref{eq:ivan}.
\eproo

\mgap

\blemm \lab{lem:dina}
Consider the sequences evolved by \emph{Algorithm \ref{alg:main}} and let $(p_k)$ and $(\what p_k)$ be as in 
\eqref{eq:def.pk}. Then, for all $k\geq 0,$
\begin{align*}
 \normq{\what p_k - p}_\gamma = (1+\alpha_k) \normq{p_k - p}_\gamma  
  - \alpha_k \normq{p_{k-1} - p}_\gamma + \alpha_k(1+\alpha_k)\normq{p_k - p_{k-1}}_\gamma \qquad \forall p \in \GG^2.
\end{align*}
\elemm
\bproo
Recall that from Lemma \ref{lm:signal}\iteref{ite:signalq} we have
\begin{align}\lab{eq:extrap4}
 \what p_k = p_k +\alpha_k(p_k - p_{k-1}), 
\end{align}
 which is clearly equivalent to
\begin{align*}
p_k - p = \dfrac{1}{1+\alpha_k} \left(\what p_k - p\right) 
 + \dfrac{\alpha_k}{1+\alpha_k}\left(p_{k-1} - p\right).
\end{align*}
Now using the well-known identity $\normq{tx+(1-t)y}_\gamma = t\normq{x}_\gamma 
+ (1-t)\normq{y}_\gamma - t(1-t)\normq{x-y}_\gamma$ with
$t=1/(1+\alpha_k)$, $x=\what p_k - p$ and $y = p_{k-1} - p$, we find
\begin{align*}
 \normq{p_k - p}_\gamma = \dfrac{1}{1+\alpha_k}\normq{\what p_k - p}_\gamma 
 + \dfrac{\alpha_k}{1+\alpha_k}\normq{p_{k-1} - p}_\gamma - \dfrac{\alpha_k}{(1+\alpha_k)^2}\normq{\what p_k - p_{k-1}}_\gamma,
\end{align*}
which, in turn, when combined with the fact that 
$\what p_k - p_{k-1} = (1+\alpha_k)(p_k - p_{k-1})$ -- see \eqref{eq:extrap4} -- and after some simple algebraic manipulations it yields
\begin{align*} 
 \normq{\what p_k - p}_\gamma = (1+\alpha_k) \normq{p_k - p}_\gamma
  - \alpha_k \normq{p_{k-1} - p}_\gamma
 + \alpha_k(1+\alpha_k)\normq{p_k - p_{k-1}}_\gamma.
\end{align*}

\eproo

\section{Asymptotic convergence of Algorithm \ref{alg:main}}
  \lab{sec:convergence}
  
In this section, we study the asymptotic convergence of Algorithm \ref{alg:main}. The main results are Theorems
\ref{th:main} and \ref{th:main02}.

\blemm \lab{lm:ineq.hk}
Consider the sequences evolved by \emph{Algorithm \ref{alg:main}} and, for an arbitrary 
$p = (z, w)\in \ES$, define
\begin{align}\lab{eq:def.hk}
h_k = \normq{p_k - p}_\gamma\qquad \forall k\geq -1.
\end{align} 
Then $h_0=h_{-1}$ and, for all $k\geq 0$,
\begin{align*}
 h_{k+1} - h_k -\alpha_k(h_k - h_{k-1}) + \tau(1 - \tau)(1 - \sigma)^2\normq{\wtil p_k - \what p_k}_\gamma
  \leq \alpha_k(1+\alpha_k)\normq{p_k - p_{k-1}}_\gamma,
\end{align*}
i.e., $(h_k)$ satisfies the assumptions of \emph{Lemma \ref{lm:alv.att}} below, where, for all $k\geq 0$,
\begin{align}
\lab{eq:def.sk}
 s_{k+1} & := \tau(1 - \tau)(1 - \sigma)^2\normq{\wtil p_k - \what p_k}_\gamma,\\[2mm]
 \lab{eq:def.deltak}
 \delta_k & := \alpha_k(1+\alpha_k)\normq{p_k - p_{k-1}}_\gamma.
\end{align}
\elemm
\bproo
The fact that $h_0 = h_{-1}$ follows directly from the fact that $p_0 = p_{-1}$ (see step 0 in Algorithm \ref{alg:main} and the definition of $p_k$ as in \eqref{eq:def.pk}).
On the other hand, from Lemma \ref{lem:dina} and the definition of $h_k$ as in \eqref{eq:def.hk},
\begin{align*} 
 \normq{\what p_k - p}_\gamma = (1+\alpha_k)\underbrace{\normq{p_k - p}_\gamma}_{h_k} 
  - \alpha_k \underbrace{\normq{p_{k-1} - p}_\gamma}_{h_{k-1}}
 + \alpha_k(1+\alpha_k)\normq{p_k - p_{k-1}}_\gamma.
\end{align*}

The desired result now follows from the above displayed equation, Lemma \ref{lm:ineq.delta}\iteref{lm:ineq.deltaq} and the definition of
$h_k$ as in \eqref{eq:def.hk}.
\eproo

\mgap

Next we present our first result on the convergence of Algorithm \ref{alg:main} when $k\to +\infty$.

\mgap

\btheo[First result on the asymptotic convergence of Algorithm \ref{alg:main}]
 \lab{th:main}
Consider the sequences evolved by \emph{Algorithm \ref{alg:main}} and let $\emptyset \neq \ES$ be as in
\eqref{eq:exts}. Assume that
\begin{align}\lab{eq:cond.sum}
\sum_{k=0}^{\infty}\,\alpha_k\normq{p_k - p_{k-1}}_\gamma< \infty,
\end{align}
where $(p_k)$ is as in \eqref{eq:def.pk}.
Then there exists $(z_\infty, w_\infty)\in \ES$ such that
\begin{align}\lab{eq:conv.zws}
z_k\to z_\infty\;\;\mbox{and}\;\; y_k\to w_\infty.
\end{align}
Additionally, we also have
\begin{align}\lab{eq:conv.v}
v_k\to z_\infty, \;\; Lx_k\to w_\infty\;\;\mbox{and}\;\;\wtil y_k\to w_\infty. 
\end{align}
\etheo
\bproo
We start by making a few remarks.
First, from \eqref{eq:cond.sum} and the fact that 
$\alpha_k (1+\alpha_k)\leq 2\alpha_k$ (because $0\leq \alpha_k<1$), we conclude that
%
$
\sum_{k=0}^\infty\,\delta_k<\infty,
$
%
where $\delta_k$ is as in \eqref{eq:def.deltak}, which, in turn, combined with Lemmas \ref{lm:ineq.hk} and \ref{lm:alv.att} (below) gives
\begin{align}\lab{eq:geremia}
 \lim_{k\to \infty}\,h_k\;\;\mbox{exists}\;\; \mbox{and}\;\;\sum_{k=1}^\infty\,s_k<\infty,
\end{align}
where $h_k$ and $s_{k+1}$ are as in \eqref{eq:def.hk} and \eqref{eq:def.sk}, respectively. 
Using \eqref{eq:def.sk} and the second statement in \eqref{eq:geremia} we also obtain
$\normq{\widetilde p_k - \widehat p_k}_\gamma\to 0$, which in turn when combined with the definitions
of $\widetilde p_k$ and $\widehat p_k$ -- as in \eqref{eq:def.pk} --, \eqref{def:norm_g3}, \eqref{eq:woodx} and \eqref{eq:wood2} yields
\begin{align}\lab{eq:geremia02}
 v_k - \widehat z_k \to 0,\;\; Lx_k - \what y_k\to 0,\;\;\; \widetilde y_k - L x_k\to 0\;\;\mbox{and}\;\; \varepsilon_k\to 0.
\end{align}

\noindent
Second, from \eqref{eq:cond.sum} and the fact that $\alpha_k^2\leq \alpha_k$, we obtain
\[
 \lim_{k\to \infty}\,\alpha_k\norm{z_k - z_{k-1}} = \lim_{k\to \infty}\,\alpha_k\norm{y_k - y_{k-1}}=0,
\]
which, in turn, when combined with the definitions of
$\what z_k$ and $\what y_k$ as in \eqref{eq:extrap} and \eqref{eq:extrap3} yields
\begin{align}\lab{eq:limits6}
\what z_k - z_k\to 0\;\;\mbox{and}\;\;\what y_k - y_k\to 0.
\end{align}

Now, let $p_k = (z_k, y_k)$ be as in \eqref{eq:def.pk}. 
Note that using the first statement in \eqref{eq:geremia}, the definition of $h_k$ as in \eqref{eq:def.hk} and Lemma \ref{lm:opial} below, it follows that to prove the  convergence of $(p_k)$ to some element in $\ES$ -- and hence the statement in \eqref{eq:conv.zws} -- it suffices to show that every cluster point of $(p_k)$ belongs to $\ES$. 
To this end, let $p_\infty = (z_\infty, y_\infty)\in \GG^2$ be a cluster point of $(p_k)$ (we know from \eqref{eq:geremia} and \eqref{eq:def.hk} that $(p_k)$ is bounded), i.e., let $z_\infty$ and $y_\infty$ be cluster points of
$(z_k)$ and $(y_k)$, respectively. 
Then let also $(k_j)$ be an increasing sequence of indexes such that 
\begin{align}\lab{eq:limits}
z_{k_j}\to z_\infty\;\;\mbox{and}\;\;y_{k_j}\to y_\infty.
\end{align}
In view of \eqref{eq:limits6} and \eqref{eq:limits}, we have
\begin{align}\lab{eq:limits2}
\what z_{k_j}\to z_\infty\;\;\mbox{and}\;\; \what y_{k_j}\to y_\infty,
\end{align}
which, in particular, when combined with \eqref{eq:geremia02} gives
\begin{align}\lab{eq:limits8}
 v_{k_j}\to z_\infty,\quad Lx_{k_j}\to y_\infty,\quad \widetilde y_{k_j}\to y_\infty\quad \mbox{and}\quad \varepsilon_{k_j}\to 0.
\end{align}

From \eqref{def:zkk}, the second statement in \eqref{eq:geremia02} (with $k = k_j$) and the first statement in \eqref{eq:limits2} we also obtain $z'_{k_j}\to z_\infty$,  which combined with Lemma \ref{lm:signal}\iteref{ite:signaln} (with $k=k_j$), the fact that
the graph of $\partial(f^*\circ -L^*)$ is closed and the second statement in \eqref{eq:limits8} yields
$-y_\infty\in \partial(f^*\circ -L^*)(z_\infty)$. 
As a consequence, according to the definition of $\mathcal{S}$ as in \eqref{eq:exts}, to prove that 
$(z_\infty, y_\infty)\in \mathcal{S}$, it remains to verify that $y_\infty\in \partial g^*(z_\infty)$. 
To this end, recall first that from \eqref{eq:wood} (with $k = k_j$) we know that 
$v_{k_j}\in \partial_{\varepsilon_{k_j}} g(\widetilde y_{k_j})$, which is equivalent to 
$\widetilde y_{k_j}\in \partial_{\varepsilon_{k_j}} g^*(v_{k_j})$.
Combining the latter inclusion with the first, third and fourth statements in \eqref{eq:limits8} as well as with
the closedness of the graph of the $\varepsilon$-subdifferential of $g^*$, we obtain the desired result, namely $y_\infty\in \partial g^*(z_\infty)$.

Altogether, we have proved that every cluster point of $(p_k)$ belongs to $\mathcal{S}$ and so, as we explained above, it guarantees that $(p_k)$ converges to some element in $\mathcal{S}$, i.e., here we finish the proof of \eqref{eq:conv.zws}.

Finally, the proof of \eqref{eq:conv.v} follows trivially from \eqref{eq:conv.zws}, \eqref{eq:geremia02} and \eqref{eq:limits6}.
\eproo

\mgap

\noindent
{\bf Remark:}
As we discussed in the Introduction (following Assumption \ref{ass:blanket}), under standard regularity conditions on \eqref{eq:prob}, the result on $(z_\infty, w_\infty)$ as in Theorem \ref{th:main}, gives that 
there exists $x_\infty\in \HH$ such that $x_\infty\in \partial f^*(-L^*z_\infty)$, $w_\infty = Lx_\infty$ and 
$x_\infty$ and $z_\infty$ are solutions of \eqref{eq:mono} and \eqref{eq:d_mono}, respectively. 
Moreover, the second statements in \eqref{eq:conv.zws} and \eqref{eq:conv.v} give, in particular, that 
$Lx_k - y_k \to 0$.

\mgap

We will consider the following two sufficient conditions on the sequences $(\alpha_k)$ and/or $(p_k)$ to ensure
\eqref{eq:cond.sum} holds -- see \eqref{def:inner_g3} and the definition of $p_k$ as in \eqref{eq:def.pk} -- :  

\mgap

\noindent
{\bf Assumption A:} for some $0<\theta<1$ and $k_0\geq 1$,
\begin{align}\label{eq:valentin2}
&\alpha_k\leq \min\left\{\alpha, \dfrac{\theta^k}{\gamma^{-1}\normq{z_k - z_{k-1}} + \gamma \normq{y_k - y_{k-1}}}\right\},\qquad \forall k\geq k_0;
\end{align}
here we adopt the convention $1/0 = \infty$.

\mgap
\mgap

\noindent
{\bf Assumption B:}
$(\alpha,\sigma,\tau)\in  [0,1)\times [0,1)\times  (0,1)$ and the sequence $(\alpha_k)$ satisfy
\begin{align}
 \label{eq:alpha}
 0\leq \alpha_k\leq \alpha_{k+1}\leq \alpha<\beta<1\qquad \forall k\geq 0,
\end{align}
where 
\begin{align}\lab{def:beta}
\beta := \dfrac{2\eta}{1 + 2\eta + \sqrt{1 + 8\eta}}
\end{align}
and
\begin{align} \lab{def:eta}
\eta := \dfrac{(1 - \tau)(1 - \sigma)^2}{4\tau}.
\end{align}

\mgap

\blemm \lab{lem:omicron}
Under the \emph{Assumption B} on \emph{Algorithm \ref{alg:main}}, define the quadratic real function
\begin{align} \label{eq:def.q}
 q(t):=(\eta - 1) t^2 - (1 + 2\eta)t + \eta \qquad \forall t\in \R.
\end{align}
Then, $q(\alpha)>0$ and, for every $p = (z, w)\in \ES$,
\begin{align}\label{eq:sum}
 \sum_{j=0}^k\norm{p_j - p_{j-1}}_\gamma^2\leq
\dfrac{2\,\norm{p_0 - p}_\gamma^2}{(1-\alpha)q(\alpha)} \qquad \forall k\geq 1.
\end{align}
\elemm
\bproo
Note first that combining Lemmas \ref{lm:ineq.hk} and \ref{lm:ineq.delta}\iteref{lm:ineq.deltaq} (second inequality) and \eqref{def:eta} we obtain, for all $k\geq 0$,
\begin{align} \lab{eq:doug02}
 h_{k+1} - h_k -\alpha_k(h_k - h_{k-1}) + \eta \normq{p_{k+1} - \what p_k}_\gamma
  \leq \alpha_k(1+\alpha_k)\normq{p_k - p_{k-1}}_\gamma,
\end{align}
where $h_k$ is as in \eqref{eq:def.hk}.
On the other hand, using Lemma \ref{lm:signal}\iteref{ite:signalq}, the Cauchy-Schwarz inequality and the Young inequality
$2ab\leq a^2+b^2$ with $a:=\norm{p_{k+1} - p_k}_\gamma$ and $b:=\norm{p_k - p_{k-1}}_\gamma$
we find
\begin{align} \lab{eq:doug}
 \nonumber
 \norm{p_{k+1} - \what p_k}_\gamma^2
   \nonumber
    &=\norm{p_{k+1} - p_k}_\gamma^2 + \alpha_k^2\norm{p_k - p_{k-1}}_\gamma^2 - 
        2\alpha_k\inner{p_{k+1} - p_k}{p_k - p_{k-1}}_\gamma\\
  \nonumber
    & \geq \norm{p_{k+1} - p_k}_\gamma^2 + \alpha_k^2\norm{p_k - p_{k-1}}_\gamma^2 - \alpha_k 
            \left(2\norm{p_{k+1} - p_k}_\gamma\norm{p_k - p_{k-1}}_\gamma\right)\\
    & \geq (1 - \alpha_k) \norm{p_{k+1} - p_k}_\gamma^2 - \alpha_k(1 - \alpha_k)\norm{p_k - p_{k-1}}_\gamma^2.
\end{align}
Using \eqref{eq:doug02}, \eqref{eq:doug} and some simple algebraic manipulations we find
\begin{align} \lab{eq:305a}
 h_{k+1} - h_k -\alpha_k(h_k - h_{k-1}) + \eta(1 - \alpha_k) \normq{p_{k+1} - p_k}_\gamma
  \leq \gamma_k\normq{p_k - p_{k-1}}_\gamma,
\end{align}
where, for all $k\geq 0$,
\begin{align}\lab{def:gamma}
\gamma_k := (1 - \eta)\alpha_k^2 + (1 + \eta)\alpha_k. 
\end{align}
Define
\begin{align} \label{eq:305}
\mu_0 := (1-\alpha_0) h_0\geq 0\;\;\mbox{and}\;\;
\mu_k := h_k - \alpha_{k-1}h_{k-1} + \gamma_{k} \norm{p_k - p_{k-1}}_\gamma^2\quad \forall k\geq 1,
\end{align}
where $h_k$ is as in \eqref{eq:def.hk}.
 Using \eqref{eq:def.q}, the assumption that $(\alpha_k)$ is nondecreasing -- see \eqref{eq:alpha} -- and \eqref{eq:305a}--\eqref{eq:305} we obtain, for all $k\geq 1$,
\begin{align}
 \label{eq:307}
 \nonumber
 \mu_k-\mu_{k-1}
\nonumber
  &\leq
\left[h_k - h_{k-1} - \alpha_{k-1}(h_{k-1} - h_{k-2})-
\gamma_{k-1}\norm{p_{k-1} - p_{k-2}}_\gamma^2\right]+\gamma_{k}\norm{p_k - p_{k-1}}_\gamma^2\\
 \nonumber
 &\leq \left[\gamma_{k} - \eta (1-\alpha_{k})\right]\norm{p_k - p_{k-1}}_\gamma^2\\
\nonumber
&=-\left[(\eta-1)\alpha_{k}^2-(1+2\eta)\alpha_{k}+\eta\right]\norm{p_k - p_{k-1}}_\gamma^2\\
&=-q(\alpha_k)\norm{p_k - p_{k-1}}_\gamma^2.
\end{align}
Note now that $0< \beta < 1$ as in \eqref{def:beta} is either the smallest or the largest root of
the quadratic function $q(\cdot)$. Hence, from \eqref{eq:alpha}, for all $k\geq 0$,
\begin{align*}
 q(\alpha_{k})\geq q(\alpha)>q(\beta) = 0.
\end{align*}
The above inequalities combined with \eqref{eq:307} yield
\begin{align}
  \label{eq:ineq.mu}
 \norm{p_k - p_{k-1}}_\gamma^2\leq \dfrac{1}{q(\alpha)}(\mu_{k-1}-\mu_{k}),\quad \forall k\geq 1,
\end{align}
which combined with \eqref{eq:alpha} and the definition of $\mu_k$ as in \eqref{eq:305} gives
\begin{align}
  \label{eq:sum.q}
 \nonumber
 \sum_{j=0}^k\,\norm{p_j - p_{j-1}}_\gamma^2&\leq \dfrac{1}{q(\alpha)}(\mu_0-\mu_k),\\
             &\leq \dfrac{1}{q(\alpha)}(\mu_0+\alpha h_{k-1}) \quad \forall k\geq 1.
\end{align}
Note now that using \eqref{eq:alpha}, \eqref{eq:305} and \eqref{eq:ineq.mu} we also find
\begin{align*}
  \mu_0\geq \ldots \geq \mu_{k} = &h_{k} - \alpha_{k-1}h_{k-1} + \gamma_{k}\|p_{k} - p_{k-1}\|_\gamma^{2} \\
  \geq&  h_{k} - \alpha h_{k-1},\quad \forall k\geq 1,
\end{align*}
and so
\begin{equation}
 \label{H13}
    h_{k}\leq \alpha^{k} h_{0}+\frac{\mu_{0}}{1-\alpha}\leq h_{0} + \dfrac{\mu_{0}}{1-\alpha} \qquad \forall k\geq 0.
\end{equation}
Hence, \eqref{eq:sum} follows directly from \eqref{eq:sum.q}, \eqref{H13}, the definition of $\mu_0$ as in \eqref{eq:305} and
the definition of $h_0$ as in \eqref{eq:def.hk}.
\eproo

\mgap

\btheo[Second result on the asymptotic convergence of Algorithm \ref{alg:main}] \lab{th:main02}
Under the assumptions \emph{A} or \emph{B} on the sequence $(\alpha_k)$, all the conclusions of
\emph{Theorem \ref{th:main}} hold true. 
\etheo
\bproo
The proof follows form Theorem \ref{th:main} (see \eqref{eq:cond.sum}), Assumptions A and B above and 
Lemma \ref{lem:omicron}.
\eproo

\section{Global convergence rates of Algorithm \ref{alg:main}}
  \lab{sec:global}
  
In this section, we study global convergence rates for Algorithm \ref{alg:main}. We obtain (global) pointwise
$O(1/\sqrt{k})$ and ergodic $O(1/k)$ rates for residuals; see Theorems \ref{the:point} and \ref{the:ergodic} below.

\mgap

\blemm \lab{lem:sat}
Consider the sequences evolved by \emph{Algorithm \ref{alg:main}} and assume that
\begin{center}
\emph{Assumption B} holds.
\end{center}
Let $(p_k)$, $(\wtil p_k)$ and 
$(\what p_k)$ be as in \eqref{eq:def.pk} and let also $q(\cdot)$ be as in \eqref{eq:def.q}. 
Then, for all $p = (z, w)\in \ES$,
\begin{align*}
\normq{p_k - p}_\gamma + \tau(1 - \tau)(1 - \sigma)^2\sum_{j=0}^k\,\normq{\wtil p_j - \what p_j}_\gamma 
\leq \left(1 + \dfrac{2\alpha(1 + \alpha)}{(1-\alpha)^2q(\alpha)}\right) \normq{p_0 - p}_\gamma\,.
\end{align*}
\elemm
\bproo
From Lemmas \ref{lm:ineq.hk} and \ref{lm:alv.att}(a) (below),
\begin{align*}
 h_k + \sum_{j=1}^k\,s_j\leq h_0 + \dfrac{1}{1 - \alpha} \sum_{j=0}^{k-1}\,\delta_j,
\end{align*}
where $h_k$, $s_k$ and $\delta_k$ are as in \eqref{eq:def.hk}, \eqref{eq:def.sk} and \eqref{eq:def.deltak}, respectively.
Then, in view of \eqref{eq:sum},
\begin{align*}
\nonumber
\normq{p_k - p}_\gamma + \tau(1 - \tau)(1 - \sigma)^2\sum_{j=0}^k\,\normq{\wtil p_j - \what p_j}_\gamma&\leq
\normq{p_0 - p}_\gamma + \dfrac{1}{1 - \alpha} \sum_{j=0}^{k-1}\,\alpha_j(1+\alpha_j)\normq{p_j - p_{j-1}}_\gamma\\[2mm]
& \leq \left(1 + \dfrac{2\alpha(1 + \alpha)}{(1-\alpha)^2q(\alpha)}\right)\normq{p_0 - p}_\gamma.
\end{align*}
\eproo

\mgap

\btheo[Pointwise global convergence rates of Algorithm \ref{alg:main}]
\lab{the:point}
Consider the sequences evolved by \emph{Algorithm \ref{alg:main}} and assume that
\begin{center}
\emph{Assumption B} holds.
\end{center}
Let $(z'_k)$ be as in \eqref{def:zkk} and let
$d_0$ denote the distance of $p_0 = (z_0, y_0)$ to the solution set $\ES$ as in \eqref{eq:exts}.
Then, for every $k\geq 0$, there exists $0\leq i\leq k$ such that
\begin{align} \lab{eq:point}
\begin{cases}
&-Lx_i \in \partial \left(f^*\circ -L^*\right)(z'_i),\qquad \wtil y_i \in \partial_{\varepsilon_i}g^*(v_i),\\[4mm]
& \gamma\normq{L x_i - \wtil y_i} + \dfrac{1}{\gamma}\normq{z'_i - v_i} \leq \dfrac{2C d_0^2}{k},\\[4mm]
& \varepsilon_i \leq \dfrac{\sigma^2 C d_0^2}{2k},
\end{cases}
\end{align}
where
\begin{align}\lab{def:C}
C := \dfrac{1}{\tau(1 - \tau)(1 - \sigma)^2}\left(1 + \dfrac{2\alpha(1 + \alpha)}{(1-\alpha)^2q(\alpha)}\right).
\end{align}
\etheo
\bproo
Let $p^* = (z^*, w^*)\in \ES$ be such that $d_0 = \norm{p_0 - p^*}_\gamma$.
From Lemma \ref{lem:sat} (with $p = p^*$) and the definition of $C>0$ as in \eqref{def:C},
\begin{align} \lab{eq:air}
\sum_{j=0}^k\,\normq{\wtil p_j - \what p_j}_\gamma \leq C d_0^2.
\end{align}
From \eqref{def:norm_g3} and Lemmas \ref{lm:signal}\iteref{ite:signalo} and 
\ref{lem:psigma}\iteref{ite:psigman},
\begin{align}\lab{eq:pal2}
\dfrac{\gamma}{2}\normq{L x_k - \wtil y_k} + \dfrac{1}{2\gamma}\normq{z'_k - v_k} 
= \dfrac{1}{2}\normq{\breve{p}_k - \what p_k}_\gamma \leq \normq{\wtil p_k - \what p_k}_\gamma.
\end{align}
Due to \eqref{def:norm_g3}, \eqref{eq:wood2} and the definitions of $\wtil p_k$ and $\what p_k$ as in \eqref{eq:def.pk} we also find
\begin{align}\lab{eq:dant3}
\dfrac{2\varepsilon_k}{\sigma^2} &\leq \gamma\normq{Lx_k - \what y_k} + \dfrac{1}{\gamma}\normq{v_k - \what z_k}
= \normq{\wtil p_k - \what p_k}_\gamma.
\end{align}
Hence, from \eqref{eq:air} -- \eqref{eq:dant3},
\begin{align}\lab{eq:madrid}
\sum_{j=0}^k\,\Delta_j \leq C d_0^2,
\end{align}
where
\begin{align} \lab{def:delta}
\Delta_j := \max\left\{\dfrac{\gamma}{2}\normq{L x_j - \wtil y_j} + \dfrac{1}{2\gamma}\normq{z'_j - v_j},
\dfrac{2\varepsilon_j}{\sigma^2} 
\right\},\qquad j = 0,\dots, k.
\end{align}
The two inequalities in \eqref{eq:point} follow by choosing $i\in \{0,\dots, k\}$ such that
$\Delta_i\leq \Delta_j$ for all $j=0, \dots, k$ and using \eqref{eq:madrid} and the definition of $\Delta_i$ as in \eqref{def:delta}.
To finish the proof of the theorem, note that the inclusions in \eqref{eq:point} follow directly from \eqref{eq:wood} (combined with the fact that $(\partial_{\varepsilon_k} g)^{-1} = \partial_{\varepsilon_k}g^*$) and Lemma \ref{lm:signal}\iteref{ite:signaln}.
\eproo

\mgap

For the sequences generated by Algorithm \ref{alg:main} and $(z'_k)$ as in \eqref{def:zkk}, define the ergodic means
\begin{align}
\lab{def:erg01}
\begin{aligned}
& x^a_k := \dfrac{1}{k+1}\sum_{j=0}^k\,x_j,\quad \wtil y^a_k := \dfrac{1}{k+1}\sum_{j=0}^k\,\wtil y_j, \\[2mm]
& z'^a_k := \dfrac{1}{k+1}\sum_{j=0}^k\,z'_j,\quad v^a_k := \dfrac{1}{k+1}\sum_{j=0}^k\,v_j.
\end{aligned}
\end{align}

Define also, for all $k\geq 0$,
\begin{align} \lab{def:erg02}
\begin{aligned}
&\delta^a_k := \dfrac{1}{k+1}\sum_{j=0}^k\,\inner{z'_j}{L(x^a_k - x_j)},\\[2mm]
&\varepsilon^a_k := \dfrac{1}{k+1}\sum_{j=0}^k\,\left[\varepsilon_j + \inner{\wtil y_j }{v_j - v^a_k}\right].
\end{aligned}
\end{align}

\mgap
\mgap

\blemm \lab{lem:erg}
Consider the sequences evolved by \emph{Algorithm \ref{alg:main}}, let $(x^a_k)$, $(\wtil y^a_k)$, $(z'^a_k)$ and $(v^a_k)$ be as in \eqref{def:erg01} and let $(\delta^a_k)$ and $(\varepsilon^a_k)$ be as in \eqref{def:erg02}. Let also
$(\breve{p}_k)$, $(\what p_k)$ and $(\wtil p_k)$ be as in \eqref{eq:def.pk}.
For all $k\geq 0$,
\benum
\item \lab{ergm}
$\delta^a_k, \varepsilon^a_k\geq 0$\; and\; 
$-Lx^a_k \in \partial_{\delta^a_k} \left(f^*\circ -L^*\right)(z'^a_k)$,\; $\wtil y^a_k\in \partial_{\varepsilon^a_k}g^*(v^a_k)$.
\item \lab{ergn}
$\delta^a_k + \varepsilon^a_k = \dfrac{1}{k+1}\sum_{j=0}^k\,\left[\varepsilon_j + \inner{\breve{p}_j - \what p_j}{\wtil p^a_k - \wtil p_j}_\gamma\right]$,
where 
\begin{align}\lab{eq:ptil}
\wtil p^a_k := \dfrac{1}{k+1}\sum_{j=0}^k\,\wtil p_j = (v^a_k, L x^a_k).
\end{align}
\eenum
\elemm
\bproo
\iteref{ergm} The desired result follows from \cite[Theorem 2.3]{mon.sva-hpe.siam10} and the inclusions in \eqref{eq:wood}
and in Lemma \ref{lm:signal}\iteref{ite:signaln}.
\mgap

\iteref{ergn} In view of Lemma \ref{lm:signal}\iteref{ite:signalo}, for $j=0,\dots, k$, we have
$ \breve{p}_j - \widehat p_j = \left(\gamma(L x_j - \widetilde y_j), \frac{1}{\gamma}\left(z'_j - v_j\right)\right)$
and so by using the definition of $(\wtil p_j)$ and \eqref{eq:ptil} we get
\begin{align*}
\sum_{j=0}^k\, \inner{\breve{p}_j - \what p_j}{\wtil p^a_k - \wtil p_j}_\gamma &=
\sum_{j=0}^k\, \left[\dfrac{1}{\gamma}\inner{\gamma(Lx_j - \wtil y_j)}{v^a_k - v_j}
+ \gamma\inner{\dfrac{1}{\gamma}(z'_j - v_j)}{Lx^a_k - Lx_j}\right]\\[2mm]
& = \sum_{j=0}^k\,\left[\inner{Lx_j - \wtil y_j}{v^a_k - v_j}
+ \inner{z'_j - v_j}{L(x^a_k - x_j)}\right]\\[2mm]
& = \sum_{j=0}^k\,\left[\inner{\wtil y_j }{v_j - v^a_k} + \inner{L x_j}{v^a_k} + \inner{z'_j}{L(x^a_k - x_j)} - \inner{v_j}{Lx^a_k}\right]\\[2mm]
& = \sum_{j=0}^k\,\left[\inner{\wtil y_j }{v_j - v^a_k} + \inner{z'_j}{L(x^a_k - x_j)}\right].
\end{align*}
The desired result now follows by adding the two equations in \eqref{def:erg02} and using the latter identity.
\eproo

\mgap
\mgap

\btheo[Ergodic global convergence rates of Algorithm \ref{alg:main}]
  \lab{the:ergodic}
Consider the sequences evolved by \emph{Algorithm \ref{alg:main}}, let $(x^a_k)$, $(\wtil y^a_k)$, $(z'^a_k)$ and $(v^a_k)$ be as in \eqref{def:erg01} and let $(\delta^a_k)$ and $(\varepsilon^a_k)$ be as in \eqref{def:erg02}.
Let also $d_0$ denote the distance of $p_0 = (z_0, y_0)$ to the solution set $\ES$ as in \eqref{eq:exts}.
Assume that $\alpha_k \equiv \alpha$ and that
\begin{center}
\emph{Assumption B} holds.
\end{center}
Then, for all $k\geq 0$, $\delta^a_k, \varepsilon^a_k\geq 0$\; and\;
\begin{align} \lab{eq:ergodic}
\begin{cases}
&-Lx^a_k \in \partial_{\delta^a_k} \left(f^*\circ -L^*\right)(z'^a_k),\qquad \wtil y^a_k\in \partial_{\varepsilon^a_k}g^*(v^a_k),\\[4mm]
& \gamma\normq{L x^a_k - \wtil y^a_k}_\gamma + \dfrac{1}{\gamma}\normq{z'^a_k - v^a_k}_\gamma \leq 
\dfrac{D^2 d_0^2}{k^2},\\[4mm]
& \delta^a_k + \varepsilon^a_k \leq \dfrac{1}{k}\left(\dfrac{\alpha(1 + \alpha)}{\tau(1-\alpha)q(\alpha)}
+ D (1 + 2\sqrt{3})\sqrt{C}\right)d_0^2,
\end{cases}
\end{align}
where $C > 0$ is as in \eqref{def:C} and
\begin{align}\lab{def:D}
D := \dfrac{1 + \alpha}{\tau}
\left(1 + \sqrt{1 + \dfrac{2\alpha(1 + \alpha)}{(1-\alpha)^2q(\alpha)}}\,\right).
\end{align}
\etheo
\bproo
Note first that the inclusions in \eqref{eq:ergodic} follow from Lemma \ref{lem:erg}\iteref{ergm}.
Now let $p^* = (z^*, w^*)\in \ES$ be such that $d_0 = \norm{p_0 - p^*}_\gamma$.
Using Lemma \ref{lm:signal}[\emph{\ref{ite:signalm}} and \emph{\ref{ite:signalq}}] and the assumption
$\alpha_k \equiv \alpha$ we find $\tau(\what p_k - \breve{p}_k) = p_k - p_{k+1} + \alpha (p_k - p_{k-1})$, for all $k\geq 0$,
and so (recall that $p_0 = p_{-1}$)
\begin{align}\lab{eq:hell}
\tau \left\|\sum_{j=0}^k\,(\what p_j - \breve{p}_j)\right\|_\gamma &= \norm{p_0 - p_{k+1}+ \alpha(p_k - p_0)}_\gamma \leq \norm{p_0 - p_{k+1}}_\gamma + \alpha\norm{p_k - p_0}_\gamma.
\end{align}
In view of Lemma \ref{lem:sat} (with $p = p^*$), the definitions of $d_0$ and $D>0$, and the triangle inequality,
\begin{align} \lab{eq:vini}
\nonumber
\norm{p_k - p_0}_\gamma &\leq \norm{p_k - p^*}_\gamma + \norm{p^* - p_0}_\gamma\\[2mm]
\nonumber
&\leq \left(1 + \sqrt{1 + \dfrac{2\alpha(1 + \alpha)}{(1-\alpha)^2q(\alpha)}}\,\right)d_0\\[2mm]
& = \dfrac{\tau D d_0}{1 + \alpha}\qquad \forall k\geq 0,
\end{align}
which combined with \eqref{def:D} and \eqref{eq:hell} yields
\begin{align} \lab{eq:emanuel}
\left\|\sum_{j=0}^k\,(\what p_j - \breve{p}_j)\right\|_\gamma & \leq D d_0.
\end{align}
Recall that from Lemma \ref{lm:signal}\iteref{ite:signalo} we have
\begin{align*}
 \breve{p}_k - \widehat p_k = \left(\gamma(L x_k - \widetilde y_k), \frac{1}{\gamma}\left(z'_k - v_k\right)\right),
\end{align*}
and so from the definitions of ergodic means as in \eqref{def:erg01}, we find
\begin{align*}
\dfrac{1}{k+1}\sum_{j=0}^k\,(\breve{p}_j - \what p_j) = \left(\gamma(L x^a_k - \wtil y^a_k), \dfrac{1}{\gamma}(z'^a_k - v^a_k)\right).
\end{align*}
Hence, in view of \eqref{def:norm_g3} and \eqref{eq:emanuel},
\begin{align*}
\gamma\normq{L x^a_k - \wtil y^a_k}_\gamma + \dfrac{1}{\gamma}\normq{z'^a_k - v^a_k}_\gamma =
      \dfrac{1}{(k+1)^2}\left\|\sum_{j=0}^k\,(\breve{p}_j - \what p_j)\right\|^2_\gamma \leq \dfrac{D^2 d_0^2}{k^2},
\end{align*}
which gives the first inequality in \eqref{eq:ergodic}.  

\mgap

Now let's prove the second inequality in \eqref{eq:ergodic}. To this end, let $p = (z, w)\in \GG^2$ and first note that from Lemma \ref{lm:ineq.delta}\iteref{lm:ineq.deltan},
\begin{align}\lab{eq:silv}
\normq{p - \what p_k}_\gamma - \normq{p - \breve{p}_k}_\gamma \geq 2\left[\varepsilon_k + \inner{\breve{p}_k - 
  \what p_k}{p - \wtil p_k}_\gamma\right].
\end{align}
By Lemma \ref{lm:signal}(a) and the identity 
$\normq{(1-\tau) a + \tau b}_\gamma = (1-\tau)\normq{a}_\gamma + \tau\normq{b}_\gamma -(1-\tau)\tau\normq{a - b}_\gamma$
with $a = p - \what p_k$ and $b = p - \breve{p}_k$,
\begin{align*}
\normq{p - p_{k+1}}_\gamma = (1-\tau)\normq{p - \what p_k}_\gamma + \tau\normq{p - \breve{p}_k}_\gamma - (1-\tau)\tau\normq{\breve{p}_k - \what p_k}_\gamma.
\end{align*}
Now multiplying \eqref{eq:silv} by $\tau > 0$ and using the latter identity we get
\begin{align} \lab{eq:lagrange}
\nonumber
\normq{p-\what p_k}_\gamma - \normq{p - p_{k+1}}_\gamma &\geq 2\tau \left[\varepsilon_k + \inner{\breve{p}_k - 
  \what p_k}{p - \wtil p_k}_\gamma\right] + (1 - \tau)\tau\normq{\breve{p}_k - \what p_k}_\gamma\\[2mm]
& \geq 2\tau \left[\varepsilon_k + \inner{\breve{p}_k - 
  \what p_k}{p - \wtil p_k}_\gamma\right].
\end{align}
Note that from Lemma \ref{lem:dina} the assumption $\alpha_k \equiv \alpha$ we obtain
\begin{align} \lab{eq:euler}
 \normq{ p - \what p_k}_\gamma = (1 + \alpha) \normq{p - p_k}_\gamma 
  - \alpha \normq{p - p_{k-1}}_\gamma
 + \alpha (1 + \alpha)\normq{p_k - p_{k-1}}_\gamma.
\end{align}
Making the substitution of \eqref{eq:euler} into \eqref{eq:lagrange} and after some simple algebra, we find (now replacing the index $k\geq 0$ by $j\geq 0$),
\begin{align*}
%
\normq{p - p_j}_\gamma - \normq{p - p_{j+1}}_\gamma + \alpha (1 + \alpha)\normq{p_j - p_{j-1}}_\gamma &\geq 2\tau \left[\varepsilon_j + \inner{\breve{p}_j - \what p_j}{p - \wtil p_j}_\gamma\right]\\[2mm] 
&\hspace{1.5cm} - \alpha \left[\normq{p - p_j}_\gamma - \normq{p - p_{j -1}}_\gamma\right].
\end{align*}
Summing up the latter inequality from $j = 0,\dots, k$ and with $p = \wtil p^a_k$ --- see \eqref{eq:ptil} --, and using 
Lemma \ref{lem:erg}\iteref{ergn},
\begin{align*}
\normq{\wtil p^a_k - p_0}_\gamma - \normq{\wtil p^a_k - p_{k+1}}_\gamma 
 &+ \alpha(1 + \alpha)\sum_{j=0}^k\,\normq{p_j - p_{j-1}}_\gamma\\
&\geq 2\tau \sum_{j=0}^k\,\left[\varepsilon_j + \inner{\breve{p}_j - \what p_j}{\wtil p^a_k - \wtil p_j}_\gamma\right] 
- \alpha \left[\normq{\wtil p^a_k - p_k}_\gamma - \normq{\wtil p^a_k - p_0}_\gamma\right]\\[2mm]
& = 2\tau (k + 1)(\delta^a_k + \varepsilon^a_k) - \alpha \left[\normq{\wtil p^a_k - p_k}_\gamma - \normq{\wtil p^a_k - p_0}_\gamma\right],
\end{align*}
which combined with \eqref{eq:sum} (with $p = p^*$) and the definition of $d_0$ yields
\begin{align*}
%
2\tau (k + 1)(\delta^a_k + \varepsilon^a_k) - \dfrac{2\alpha(1 + \alpha)d_0^2}{(1-\alpha)q(\alpha)}&\leq 
 \Big[\normq{\wtil p^a_k - p_0}_\gamma - \normq{\wtil p^a_k - p_{k+1}}_\gamma\Big]\\[2mm]
&\hspace{1.5cm}+ \alpha \Big[\normq{\wtil p^a_k - p_k}_\gamma - \normq{\wtil p^a_k - p_0}_\gamma\Big].
\end{align*}
Now using the inequality $\normq{a}_\gamma - \normq{b}_\gamma\leq 2\norm{a}_\gamma\norm{a - b}_\gamma$ (in both terms in the right-hand side of the latter inequality)  and \eqref{eq:vini} we find
\begin{align} \lab{eq:clapton}
\nonumber
2\tau (k + 1)(\delta^a_k + \varepsilon^a_k) - \dfrac{2\alpha(1 + \alpha)\,d_0^2}{(1-\alpha)q(\alpha)}&\leq
  2\norm{\wtil p^a_k - p_0}_\gamma \norm{p_{k+1} - p_0}_\gamma 
     + 2\alpha\norm{\wtil p^a_k - p_k}_\gamma \norm{p_k - p_0}_\gamma\\[2mm]
\nonumber
 &\leq \dfrac{2\tau D d_0}{1 + \alpha}
\Big(\norm{\wtil p^a_k - p_0}_\gamma + \alpha\norm{\wtil p^a_k - p_k}_\gamma\Big)\\[2mm]
& \leq 2\tau D d_0 \max\{\norm{\wtil p^a_k - p_0}_\gamma, \norm{\wtil p^a_k - p_k}_\gamma\}.
\end{align}

Now define, for all $k\geq 0$,
\begin{align} \lab{def:puk}
 \breve{p}^a_k := \dfrac{1}{k + 1}\sum_{j=0}^k\,\breve{p}_j.
\end{align}
Using \eqref{eq:wood2}, \eqref{eq:def.pk}, Lemma \ref{lem:psigma}\iteref{ite:psigmam}, \eqref{eq:ptil}, \eqref{def:puk} 
as well as the convexity of $\normq{\cdot}_\gamma$ we find
\begin{align} \lab{eq:doce}
\nonumber
\normq{\wtil p^a_k - \breve{p}^a_k}_\gamma &\leq \dfrac{1}{k + 1}\sum_{j=0}^k\,\normq{\wtil p_j - \breve{p}_j}_\gamma\\[2mm]
\nonumber 
 & = \dfrac{1}{k + 1}\sum_{j=0}^k\,\Big[\dfrac{1}{\gamma}\normq{e_j} + 
      \dfrac{1}{\gamma}\normq{v_j - \what z_j}\Big]\\[2mm]
\nonumber
&\leq \dfrac{1}{k+1}\sum_{j=0}^k\,\Big[ \gamma \normq{Lx_j - \what y_j} +  \dfrac{1}{\gamma}\normq{v_j - \what z_j}\Big]\\[2mm]
\nonumber
& = \dfrac{1}{k+1}\sum_{j=0}^k\,\normq{\wtil p_j - \what p_j}_\gamma\\[2mm]
\nonumber
& \leq \dfrac{C d_0^2}{k + 1}\\[2mm]
&\leq Cd_0^2,
\end{align}
where we used Lemma \ref{lem:sat} (with $p = p^*$) and the definition of $C>0$ as in \eqref{def:C}.
Similarly, using \eqref{def:puk}, the triangle inequality and the well-known inequality 
$\normq{a + b}_\gamma\leq 2\left(\normq{a}_\gamma + \normq{b}_\gamma\right)$, we get, for all $\ell \geq 1$,
\begin{align} \lab{eq:coff}
\nonumber
\normq{p_\ell - \breve{p}^a_k}_\gamma &= \left\|\dfrac{1}{k + 1}\sum_{j = 0}^k\,(p_\ell - \breve{p}_j)\right\|^2_\gamma\\[2mm]
\nonumber
&\leq \dfrac{1}{k + 1}\sum_{j=0}^k\,\normq{p_\ell - \breve{p}_j}_\gamma\\[2mm]
\nonumber
&\leq \dfrac{2}{k + 1}\sum_{j=0}^k\,\Big(\normq{p_\ell - p_{j+1}}_\gamma + \normq{p_{j+1} - \breve{p}_j}_\gamma\Big)\\[2mm]
& = \dfrac{2}{k + 1}\sum_{j=0}^k\,\normq{p_\ell - p_{j+1}}_\gamma + \dfrac{2}{k + 1}\sum_{j=0}^k\,\normq{p_{j+1} - \breve{p}_j}_\gamma.
\end{align}
Using (again) Lemma \ref{lem:sat} (with $p = p^*$) as well as the definition of $C>0$ as in \eqref{def:C}, we find, for all $\ell, j\geq 0$,
\begin{align} \lab{eq:coff02}
\normq{p_\ell - p_{j+1}}_\gamma \leq 2\Big(\normq{p_\ell - p^*}_\gamma + \normq{p_{j+1} - p^*}_\gamma\Big) 
\leq 2\tau(1 - \tau)(1 - \sigma)^2 C d_0^2.
\end{align}
Recall also that from Lemma \ref{lm:signal}\iteref{ite:signaln} we have $p_{j+1} - \breve{p}_j = (1 - \tau)(\what p_j - \breve{p_j})
$ ($0\leq j\leq k$) and so from Lemma \ref{lem:psigma}\iteref{ite:psigmam}, for all $j\geq 0$,
\begin{align*}
\norm{p_{j+1} - \breve{p}_j}_\gamma = (1- \tau)\norm{\what p_j - \breve{p_j}}_\gamma \leq 2(1-\tau)\norm{\wtil p_j - \what p_j}_\gamma,
\end{align*}
which then yields (by Lemma \ref{lem:sat} (with $p = p^*$) and the definition of $C > 0$)
\begin{align} \lab{eq:coff03}
\sum_{j=0}^k\,\normq{p_{j+1} - \breve{p}_j}_\gamma &= 4(1 - \tau)^2\sum_{j=0}^k\,\normq{\wtil p_j - \what p_j}_\gamma \leq 4(1 - \tau)^2Cd_0^2.
\end{align}
Putting it all together, from \eqref{eq:coff} -- \eqref{eq:coff03}  we obtain, for all $\ell \geq 0$,
\begin{align}\lab{eq:coff04}
\nonumber
\normq{p_\ell - \breve{p}^a_k}_\gamma &\leq 4 \tau(1 - \tau)(1 - \sigma)^2 C d_0^2
+ \dfrac{8}{k+1}(1 - \tau)^2Cd_0^2\\[2mm]
\nonumber
& = 4(1 - \tau)\left( \tau(1 - \sigma^2)
+ \dfrac{2}{k + 1}(1 - \tau)\right)Cd_0^2\\[2mm]
& \leq 12 C d_0^2.
\end{align}
Using now the triangle inequality, \eqref{eq:doce} and \eqref{eq:coff04} we find, for all $\ell \geq 0$,
\begin{align}\lab{eq:coff05}
\norm{p_\ell - \wtil p^a_k}_\gamma&\leq \norm{p_\ell - \breve{p}^a_k}_\gamma + \norm{\breve{p}^a_k - \wtil p^a_k}_\gamma 
\leq (1 + 2\sqrt{3})\sqrt{C}d_0.
\end{align}
Finally, using \eqref{eq:clapton} and \eqref{eq:coff05} (with $\ell = 0$ and $\ell = k$),
\begin{align*}
2\tau (k + 1)(\delta^a_k + \varepsilon^a_k) - \dfrac{2\alpha(1 + \alpha)\,d_0^2}{(1-\alpha)q(\alpha)}
 \leq 2\tau D (1 + 2\sqrt{3})\sqrt{C} d_0^2
\end{align*}
so that
\begin{align*}
\delta^a_k + \varepsilon^a_k\leq \dfrac{1}{k}\left(\dfrac{\alpha(1 + \alpha)}{\tau(1-\alpha)q(\alpha)}
+ D (1 + 2\sqrt{3})\sqrt{C}\right)d_0^2.
\end{align*}
\eproo

\section{Numerical Experiments}
  \label{sec:num}
This section presents some numerical experiments on the LASSO problem, which is an instance of the minimization problem \eqref{eq:prob}.
We compared Algorithm \ref{alg:main} from this paper with and without inertial effects; they are called \textit{Inexact ADMM} and \textit{Inexact inertial ADMM}, respectively. 
We implemented both algorithms in Matlab R2021a and, for both algorithms and all problem classes, used the same stopping criterion, namely
\begin{align}\label{eq:valentin}
 \mbox{dist}_{\infty}\big(0,\partial f(x_k)+ \partial g(x_k)\big)\leq \varepsilon,
\end{align}
where $\mbox{dist}_{\infty}(0,S):=\inf\{\|s\|_\infty\,|\,s\in S\}$ and $\varepsilon$ is a tolerance parameter set to $10^{-6}$.

The inertial parameter $\alpha_k$ (as in step 1 of Algorithm \ref{alg:main}) is updated 
according to the rule \eqref{eq:valentin2} with $\theta=0.99$ and $k_0=1$. More precisely, we choose
$\alpha_k$ as
\begin{align*}
    \alpha_k = \min\left\{\alpha, \dfrac{\theta^k}{\gamma^{-1}\normq{z_k - z_{k-1}} + \gamma \normq{y_k - y_{k-1}}}\right\},\quad \forall k\geq 1,
\end{align*}
where $0\leq \alpha<1$. The source codes are available under request (marina.geremia@ifsc.edu.br).

\mgap

\noindent
{\bf The LASSO problem.} We perform numerical experiments on the LASSO problem (as already discussed in \eqref{eq:lasso}),
namely
\begin{align}\lab{eq:lassob}
\min_{x\in \R^d}\,\left\{\dfrac{1}{2}\normq{Ax - b}+\nu\norm{x}_1\right\},
\end{align}
where $A\in \R^{n\times d}$, $b\in \R^n$ and $\nu>0$.
For the data matrix $A$ and the vector $b$, we used five categories of non-artificial datasets 
(available at the UCI Machine Learning Repository, https://archive.ics.uci.edu):
\begin{itemize}
\item[] \emph{BlogFeedback}: This category consists of one standard microarray datasets that contain features extracted from a blog post from~\cite{misc_blogfeedback_304}. This problem is called \emph{blogFeedback} (with $n=60021$ and $d=280$).
 \item[] \emph{Breast Cancer Wisconsin (Prognostic)}: This category consists of one prognostic Wisconsin breast cancer database from~\cite{misc_breast_cancer_wisconsin_(prognostic)_16}, which has a dense matrix $A$. Each row represents follow-up data for one breast cancer case. This problem is called \emph{Wisconsin} (with $n=198$ and $d=33$).
 \item[] \emph{DrivFace}: This category comprises a single standard microarray dataset containing image sequences of individuals driving in real-world scenarios from ~\cite{misc_drivface_378}. This problem is called \emph{DrivFace} (with $n=606$ and $d=6400$) and has a dense matrix $A$.
 \item[] \emph{Gene expression}: This category consists of six standard cancer DNA microarray datasets from~\cite{det.buh-fin.jma04}, which have dense and wide matrices $A$, with the number of rows 
 $n \in [42, 102]$ and the number of columns $d \in [2000, 6033]$. These problems are called \emph{brain} (with $n=42$ and $d=5597$), \emph{colon} (with $n=62$ and $d=2000$), \emph{leukemia} (with $n=72$ and $d=3571$), 
 \emph{lymphoma} (with $n=62$ and $d=4026$), \emph{prostate} (with $n=102$ and $d=6033$) and 
 \emph{srbct} (with $n=63$ and $d=2308$).
  \item[] \emph{Single-Pixel camera}: This category consists of four compressed image sensing datasets from~\cite{singlepixel}, which have dense and wide matrices $A$, with $n \in \{410, 1638\}$ and $d\in \{1024, 4096\}]$. 
 These problems are called \emph{Ball64\_singlepixcam} (with $n=1638$ and $d=4096$), \emph{Logo64\_singlepixcam} (with $n=1638$ and $d=4096$), \emph{Mug32\_singlepixcam} (with $n=410$ and $d=1024$) and \emph{Mug128\_singlepixcam} (with $n=410$ and $d=1024$).
 \end{itemize}

We implemented both algorithms \textit{Inexact ADMM} and \textit{Inexact inertial ADMM} in Matlab R2021a, combined with a CG procedure to approximately solve the subproblems \eqref{eq:light}; see also the fourth remark following Definition \ref{def:app_sol}. 
As usual (see, e.g., \cite{alv.eck-admm.coap20}), we solved the (easy) subproblem \eqref{eq:borwein} by using the 
standard-soft thresholding operator.
We also set $(\alpha, \sigma, \tau, \gamma) = (0.33, 0.99, 0.999, 1)$ and
$(\sigma, \tau, \gamma) = (0.99, 0.999, 1)$ for \textit{Inexact inertial ADMM} and \textit{Inexact ADMM}, respectively.
Moreover, as in~\cite{boy.par.chu-dis.ftml11}, we set the regularization parameter $\nu$ as $0.1\norm{A^{T}b}_{\infty}$, and scaled the vector $b$ and the columns of matrix $A$ to have $\ell_2$ unit norm. 

Table~\ref{tab:outer} shows the number of outer iterations required by each algorithm on each problem instance, the cumulative total number of inner iterations required by the CG algorithm for solving \eqref{eq:light} and runtimes in seconds demanded by each algorithm to achieve the prescribed tolerance as in \eqref{eq:valentin}.
From Table~\ref{tab:outer}, we see that \textit{Inexact inertial ADMM} outperforms \textit{Inexact ADMM} on average by about $30\%$, $25\%$ and $25\%$ on ``Outer iterations'', ``Total inner iterations'' and ``Runtimes'', respectively.

\begin{table}
\caption{Comparison of performance in the LASSO problem}
\label{tab:outer}
\scalebox{0.85}{
\begin{tabular}{|l|ccc|ccc|ccc|}
\hline
 & \multicolumn{3}{c|}{Inexact ADMM} & \multicolumn{3}{c|}{Inexact inertial ADMM} & \multicolumn{3}{l|}{} \\ \hline
 & \multicolumn{1}{c|}{\begin{tabular}[c]{@{}c@{}}Outer\\ iterations\\ \textit{(outer1)}\end{tabular}} & \multicolumn{1}{c|}{\begin{tabular}[c]{@{}c@{}}Total inner\\ iterations\\ \textit{(inner1)}\end{tabular}} & \begin{tabular}[c]{@{}c@{}}Time\\ \textit{(time1)}\end{tabular} & \multicolumn{1}{c|}{\begin{tabular}[c]{@{}c@{}}Outer \\ iterations\\ \textit{(outer2)}\end{tabular}} & \multicolumn{1}{c|}{\begin{tabular}[c]{@{}c@{}}Total inner\\ iterations\\ \textit{(inner2)}\end{tabular}} & \begin{tabular}[c]{@{}c@{}}Time\\ \textit{(time2)}\end{tabular} & \multicolumn{1}{c|}{$\frac{outer 2}{outer 1}$} & \multicolumn{1}{c|}{$\frac{inner2}{inner1}$} & $\frac{time2}{time1}$ \\ \hline
Brain & \multicolumn{1}{c|}{2923} & \multicolumn{1}{c|}{14007} & 3.97 & \multicolumn{1}{c|}{2120} & \multicolumn{1}{c|}{11112} & 3.11 & \multicolumn{1}{c|}{0.7253} & \multicolumn{1}{c|}{0.7933} & 0.7834 \\ \hline
Colon & \multicolumn{1}{c|}{505} & \multicolumn{1}{c|}{2818} & 0.39 & \multicolumn{1}{c|}{347} & \multicolumn{1}{c|}{1866} & 0.27 & \multicolumn{1}{c|}{0.6871} & \multicolumn{1}{c|}{0.6622} & 0.6923 \\ \hline
Leukemia & \multicolumn{1}{c|}{764} & \multicolumn{1}{c|}{3695} & 0.83 & \multicolumn{1}{c|}{544} & \multicolumn{1}{c|}{2861} & 0.62 & \multicolumn{1}{c|}{0.7121} & \multicolumn{1}{c|}{0.7743} & 0.7469 \\ \hline
Lymphoma & \multicolumn{1}{c|}{1101} & \multicolumn{1}{c|}{6091} & 1.45 & \multicolumn{1}{c|}{862} & \multicolumn{1}{c|}{5119} & 1.13 & \multicolumn{1}{c|}{0.7829} & \multicolumn{1}{c|}{0.8404} & 0.7793 \\ \hline
Prostate & \multicolumn{1}{c|}{2006} & \multicolumn{1}{c|}{8328} & 4.13 & \multicolumn{1}{c|}{1414} & \multicolumn{1}{c|}{6561} & 3.46 & \multicolumn{1}{c|}{0.7049} & \multicolumn{1}{c|}{0.7878} & 0.8378 \\ \hline
Srbct & \multicolumn{1}{c|}{511} & \multicolumn{1}{c|}{3554} & 0.55 & \multicolumn{1}{c|}{346} & \multicolumn{1}{c|}{2293} & 0.38 & \multicolumn{1}{c|}{0.6771} & \multicolumn{1}{c|}{0.6452} & 0.6909 \\ \hline
Ball64 & \multicolumn{1}{c|}{313} & \multicolumn{1}{c|}{490} & 8.77 & \multicolumn{1}{c|}{216} & \multicolumn{1}{c|}{325} & 5.93 & \multicolumn{1}{c|}{0.6901} & \multicolumn{1}{c|}{0.6633} & 0.6762 \\ \hline
Logo64 & \multicolumn{1}{c|}{316} & \multicolumn{1}{c|}{495} & 8.47 & \multicolumn{1}{c|}{221} & \multicolumn{1}{c|}{359} & 5.98 & \multicolumn{1}{c|}{0.6994} & \multicolumn{1}{c|}{0.7253} & 0.7061 \\ \hline
Mug32 & \multicolumn{1}{c|}{134} & \multicolumn{1}{c|}{282} & 0.09 & \multicolumn{1}{c|}{88} & \multicolumn{1}{c|}{197} & 0.06 & \multicolumn{1}{c|}{0.6567} & \multicolumn{1}{c|}{0.6986} & 0.6667 \\ \hline
Mug128 & \multicolumn{1}{c|}{955} & \multicolumn{1}{c|}{1163} & 248.04 & \multicolumn{1}{c|}{822} & \multicolumn{1}{c|}{986} & 212.21 & \multicolumn{1}{c|}{0.8607} & \multicolumn{1}{c|}{0.8478} & 0.8555 \\ \hline
DrivFace & \multicolumn{1}{c|}{2682} & \multicolumn{1}{c|}{18727} & 165.47 & \multicolumn{1}{c|}{1803} & \multicolumn{1}{c|}{13284} & 115.36 & \multicolumn{1}{c|}{0.6723} & \multicolumn{1}{c|}{0.7094} & 0.6972 \\ \hline
Wisconsin & \multicolumn{1}{c|}{285} & \multicolumn{1}{c|}{605} & 0.07 & \multicolumn{1}{c|}{182} & \multicolumn{1}{c|}{456} & 0.05 & \multicolumn{1}{c|}{0.6386} & \multicolumn{1}{c|}{0.7537} & 0.7143 \\ \hline
blogFeedback & \multicolumn{1}{c|}{386} & \multicolumn{1}{c|}{1108} & 46.79 & \multicolumn{1}{c|}{317} & \multicolumn{1}{c|}{938} & 38.84 & \multicolumn{1}{c|}{0.8212} & \multicolumn{1}{c|}{0.8466} & 0.8301 \\ \hline
Geometric mean & \multicolumn{1}{c|}{662.68} & \multicolumn{1}{c|}{2187.47} & 3.21 & \multicolumn{1}{c|}{473.79} & \multicolumn{1}{c|}{1633.23} & 2.38 & \multicolumn{1}{c|}{0.7149} & \multicolumn{1}{c|}{0.7466} & 0.7414 \\ \hline
\end{tabular}}
\end{table}

\appendix
\section{Auxiliary results}
 \lab{sec:app}

The following lemma was essentially proved by Alvarez and Attouch in  \cite[Theorem 2.1]{alv.att-iner.svva01} (see also
\cite[Lemma A.4]{alv.mar-ine.svva2020}).
\begin{lemma}
 \label{lm:alv.att}
Let the sequences $(h_k)$, $(s_k)$, $(\alpha_k)$ and $(\delta_k)$ in $[0,\infty)$
and $\alpha\in \R$ be such that
$h_0=h_{-1}$, $0\leq \alpha_k\leq \alpha<1$ and
\begin{align}
  \label{eq:alv.att02}
h_{k+1}-h_k+s_{k+1}\leq \alpha_k(h_k-h_{k-1})+\delta_k\qquad \forall k\geq 0.
\end{align}
The following hold:
\begin{enumerate}
  \item [\emph{(a)}] For all $k\geq 1$,
 \begin{align}
       \label{eq:alv.att01}
       h_k+\sum_{j=1}^k\,s_j\leq
      h_0+\dfrac{1}{1-\alpha} \sum_{j=0}^{k-1}\,\delta_j.
    \end{align}
  \item [\emph{(b)}] If $\sum^{\infty}_{k=0}\delta_k <\infty$, then $\lim_{k\to \infty}\,h_{k}$ exists, i.e., the sequence $(h_k)$ converges to some element in $[0,\infty)$.
\end{enumerate}
\end{lemma}

\mgap

\begin{lemma}[Opial~\cite{opial-wea.bams67}]
 \label{lm:opial}
Let $\HH$ be a finite dimensional inner product space, let $\emptyset \neq \mathcal{S}\subset \HH$ and let $\{p_k\}$ be a sequence in $\HH$
such that every cluster point of $\{p_k\}$ belongs to $\mathcal{S}$ and $\lim_{k\to
\infty}\,\norm{p_k - p}$ exists for every $p\in \mathcal{S}$. Then $\{p_k\}$
converges to a point in $\mathcal{S}$.
\end{lemma}


%

\def\cprime{$'$} \def\cprime{$'$}

\end{document}